\documentclass[reqno]{amsart}
\usepackage{setspace}
\usepackage{graphicx}
\usepackage{amsmath}
\usepackage{amsfonts}
\usepackage{amssymb}

\begin{document}

\centerline{\huge \bf Characterization of the Two-Dimensional}

\medskip
\centerline{\huge \bf Six-Fold Lattice Tiles}

\bigskip
\centerline{\large Chuanming Zong}

\vspace{0.6cm}
\centerline{\begin{minipage}{12.5cm}
{\bf Abstract.} This paper characterizes all the convex domains which can form six-fold lattice tilings of the Euclidean plane. They are parallelograms, centrally symmetric hexagons, centrally symmetric octagons (under suitable affine linear transformations) with vertices ${\bf v}_1=(\alpha-1 , 2)$, ${\bf v}_2=(\alpha , -2)$, ${\bf v}_3=(1-\alpha , 0)$, ${\bf v}_4=(1+\alpha , -1)$, ${\bf v}_5=-{\bf v}_1$, ${\bf v}_6=-{\bf v}_2$, ${\bf v}_7=-{\bf v}_3$ and ${\bf v}_8=-{\bf v}_4$, where $0<\alpha <{1\over 6}$, centrally symmetric decagons (under suitable affine linear transformations) with ${\bf u}_1=(-1, {1\over 2})$, ${\bf u}_2=({1\over 2}, 1)$, ${\bf u}_3=({3\over 2}, 1)$, ${\bf u}_4=(2, {1\over 2})$, ${\bf u}_5=(2,0)$, ${\bf u}_6=-{\bf u}_1$, ${\bf u}_7=-{\bf u}_2$, ${\bf u}_8=-{\bf u}_3$, ${\bf u}_9=-{\bf u}_4$ and ${\bf u}_{10}=-{\bf u}_5$ as the middle points of their edges, or with ${\bf u}_1=({1\over 2}, -1)$, ${\bf u}_2=({3\over 2}, -{1\over 2})$, ${\bf u}_3=(2, 0)$, ${\bf u}_4=({3\over 2}, {1\over 2})$, ${\bf u}_5=({1\over 2}, 1)$, ${\bf u}_6=-{\bf u}_1$, ${\bf u}_7=-{\bf u}_2$, ${\bf u}_8=-{\bf u}_3$, ${\bf u}_9=-{\bf u}_4$ and ${\bf u}_{10}=-{\bf u}_5$ as the middle points of their edges.
\end{minipage}}


\vspace{0.8cm}
\noindent
{\Large\bf 1. Introduction}

\bigskip\noindent
Let $K$ be a convex body with interior ${\rm int}(K)$, boundary $\partial (K)$ and volume ${\rm vol}(K)$, and let $X$ be a discrete set, both in $\mathbb{E}^n$. We call $K+X$ a {\it translative tiling} of $\mathbb{E}^n$ and call $K$ a {\it translative tile} if $K+X=\mathbb{E}^n$ and the translates ${\rm int}(K)+{\bf x}_i$ are pairwise disjoint. In other words, if $K+X$ is both a packing and a covering in $\mathbb{E}^n$. In particular, we call $K+\Lambda$ a {\it lattice tiling} of $\mathbb{E}^n$ and call $K$ a {\it lattice tile} if $\Lambda $ is an $n$-dimensional lattice.

Let $X$ be a discrete multiset in $\mathbb{E}^n$ and let $k$ be a positive integer. We call $K+X$ a {\it $k$-fold translative tiling} of $\mathbb{E}^n$ and call $K$ a {\it $k$-fold translative tile} if every point ${\bf x}\in \mathbb{E}^n$ belongs to at least $k$ translates of $K$ in $K+X$ and every point ${\bf x}\in \mathbb{E}^n$ belongs to at most $k$ translates of ${\rm int}(K)$ in ${\rm int}(K)+X$. In other words, $K+X$ is both a $k$-fold packing and a $k$-fold covering in $\mathbb{E}^n$. In particular, we call $K+\Lambda$ a {$k$-fold lattice tiling} of $\mathbb{E}^n$ and call $K$ a {\it $k$-fold lattice tile} if $\Lambda $ is an $n$-dimensional lattice. Clearly, a {\it $k$-fold translative tile} must be a polytope. In fact, it is proved by Gravin, Robins and Shiryaev \cite{grs} that a {\it $k$-fold translative tile} must be a centrally symmetric polytope with centrally symmetric facets.

Let $P$ be an $n$-dimensional centrally symmetric convex polytope, let $\tau (P)$ denote the smallest integer $k$ such that $P$ is a $k$-fold translative tile, and let $\tau^* (P)$ denote the smallest integer $k$ such that $P$ is a $k$-fold lattice tile. For convenience, we define $\tau (P)=\infty $ if $P$ can not form translative tiling of any multiplicity. Clearly, for every convex polytope we have
$$\tau (P)\le \tau^*(P).$$
If $\sigma $ is a non-singular affine linear transformation from $\mathbb{E}^n$ to $\mathbb{E}^n$, it can be easily verified that $P+X$ is a $k$-fold tiling of $\mathbb{E}^n$ if and only if $\sigma (P)+\sigma (X)$ is a $k$-fold tiling of $\mathbb{E}^n$. Thus, both $\tau (\sigma (P))=\tau (P)$ and $\tau^*(\sigma (P))=\tau^*(P)$ hold for all convex polytopes $P$ and all non-singular affine linear transformations $\sigma$.

\medskip
Recently, Yang and Zong \cite{yz1} proved the following results: {\it Besides parallelograms and centrally symmetric hexagons, there is no other convex domain which can form a two-, three- or four-fold lattice tiling in the Euclidean plane. However, there are convex octagons and decagons which can form five-fold lattice tilings.} Afterwards, Zong \cite{zong} characterized all the two-dimensional five-fold lattice tiles. They are
parallelograms, centrally symmetric hexagons, centrally symmetric octagons (under a suitable affine linear transformation) with vertices ${\bf v}_1=(-\alpha , -{3\over 2})$, ${\bf v}_2=(1-\alpha , -{3\over 2})$, ${\bf v}_3=(1+\alpha , -{1\over 2})$, ${\bf v}_4=(1-\alpha , {1\over 2})$, ${\bf v}_5=-{\bf v}_1$, ${\bf v}_6=-{\bf v}_2$, ${\bf v}_7=-{\bf v}_3$ and ${\bf v}_8=-{\bf v}_4$, where $0<\alpha <{1\over 4},$ or with vertices
${\bf v}_1=(\beta , -2)$, ${\bf v}_2=(1+\beta , -2)$, ${\bf v}_3=(1-\beta , 0)$, ${\bf v}_4=(\beta , 1)$, ${\bf v}_5=-{\bf v}_1$, ${\bf v}_6=-{\bf v}_2$, ${\bf v}_7=-{\bf v}_3$, ${\bf v}_8=-{\bf v}_4$, where ${1\over 4}<\beta <{1\over 3}$, or centrally symmetric decagons (under a suitable affine linear transformation) with ${\bf u}_1=(0, 1)$, ${\bf u}_2=(1, 1)$, ${\bf u}_3=({3\over 2}, {1\over 2})$, ${\bf u}_4=({3\over 2}, 0)$, ${\bf u}_5=( 1,-{1\over 2})$, ${\bf u}_6=-{\bf u}_1$, ${\bf u}_7=-{\bf u}_2$, ${\bf u}_8=-{\bf u}_3$, ${\bf u}_9=-{\bf u}_4$ and ${\bf u}_{10}=-{\bf u}_5$ as the middle points of its edges.

\smallskip
This paper characterizes all the two-dimensional six-fold lattice tiles by proving the following results.

\medskip\noindent
{\bf Theorem 1.} {\it A convex domain can form a six-fold lattice tiling of the Euclidean plane if and only if it is a parallelogram, a centrally symmetric hexagon, a centrally symmetric octagon (under suitable affine linear transformations) with vertices ${\bf v}_1=(\alpha-1 , 2)$, ${\bf v}_2=(\alpha , -2)$, ${\bf v}_3=(1-\alpha , 0)$, ${\bf v}_4=(1+\alpha , -1)$, ${\bf v}_5=-{\bf v}_1$, ${\bf v}_6=-{\bf v}_2$, ${\bf v}_7=-{\bf v}_3$ and ${\bf v}_8=-{\bf v}_4$, where $0<\alpha <{1\over 6}$, a centrally symmetric decagons (under suitable affine linear transformations) with ${\bf u}_1=(-1, {1\over 2})$, ${\bf u}_2=({1\over 2}, 1)$, ${\bf u}_3=({3\over 2}, 1)$, ${\bf u}_4=(2, {1\over 2})$, ${\bf u}_5=(2,0)$, ${\bf u}_6=-{\bf u}_1$, ${\bf u}_7=-{\bf u}_2$, ${\bf u}_8=-{\bf u}_3$, ${\bf u}_9=-{\bf u}_4$ and ${\bf u}_{10}=-{\bf u}_5$ as the middle points of their edges, or
with ${\bf u}_1=({1\over 2}, -1)$, ${\bf u}_2=({3\over 2}, -{1\over 2})$, ${\bf u}_3=(2, 0)$, ${\bf u}_4=({3\over 2}, {1\over 2})$, ${\bf u}_5=({1\over 2}, 1)$, ${\bf u}_6=-{\bf u}_1$, ${\bf u}_7=-{\bf u}_2$, ${\bf u}_8=-{\bf u}_3$, ${\bf u}_9=-{\bf u}_4$ and ${\bf u}_{10}=-{\bf u}_5$ as the middle points of their edges.}

\medskip\noindent
{\bf Theorem 2.} {\it Let $Q$ denote the quadrilateral with vertices ${\bf q}_1=(0, 1)$, ${\bf q}_2=(0, {5\over 6})$, ${\bf q}_3=(-{1\over 4}, {3\over 4})$ and ${\bf q}_4=(-{1\over 3}, {5\over 6})$. A centrally symmetric convex decagon $P_{10}$ with ${\bf u}_1=(-1, {1\over 2})$, ${\bf u}_2=({1\over 2}, 1)$, ${\bf u}_3=({3\over 2}, 1)$, ${\bf u}_4=(2, {1\over 2})$, ${\bf u}_5=(2,0)$, ${\bf u}_6=-{\bf u}_1$, ${\bf u}_7=-{\bf u}_2$, ${\bf u}_8=-{\bf u}_3$, ${\bf u}_9=-{\bf u}_4$ and ${\bf u}_{10}=-{\bf u}_5$ as the middle points of its edges if and only if one of its vertices is an interior point of $Q$.

Let $Q^*$ denote the quadrilateral with vertices ${\bf q}_1=(0, {5\over 4})$, ${\bf q}_2=({1\over 6}, {7\over 6})$, ${\bf q}_3=(0, 1)$ and ${\bf q}_4=(-{1\over 6}, {7\over 6})$. A centrally symmetric convex decagon $P^*_{10}$ with ${\bf u}_1=({1\over 2}, -1)$, ${\bf u}_2=({3\over 2}, -{1\over 2})$, ${\bf u}_3=(2, 0)$, ${\bf u}_4=({3\over 2}, {1\over 2})$, ${\bf u}_5=({1\over 2}, 1)$, ${\bf u}_6=-{\bf u}_1$, ${\bf u}_7=-{\bf u}_2$, ${\bf u}_8=-{\bf u}_3$, ${\bf u}_9=-{\bf u}_4$ and ${\bf u}_{10}=-{\bf u}_5$ as the middle points of their edges if and only if one of its vertices is an interior point of $Q^*$.}

\medskip
\noindent
{\bf Remark 1.} This paper follows the key idea and the methodology of Zong \cite{zong}.

\vspace{0.6cm}
\noindent
{\Large\bf 2. Basic Results}

\bigskip\noindent
Let $P_{2m}$ denote a centrally symmetric convex $2m$-gon centered at the origin, let ${\bf v}_1$, ${\bf v}_2$, $\ldots$, ${\bf v}_{2m}$ be the $2m$ vertices of $P_{2m}$ enumerated clock-wise, and let $G_1$, $G_2$, $\ldots $, $G_{2m}$ be the $2m$ edges of $P_{2m}$, where $G_i$ has two ends ${\bf v}_i$ and ${\bf v}_{i+1}$. For convenience, we write $V=\{{\bf v}_1, {\bf v}_2, \ldots, {\bf v}_{2m}\}$ and $\Gamma=\{G_1, G_2, \ldots, G_{2m}\}$.

Assume that $P_{2m}+X$ is a $\tau (P_{2m})$-fold translative tiling of $\mathbb{E}^2$, where $X=\{{\bf x}_1, {\bf x}_2, {\bf x}_3, \ldots \}$ is a discrete multiset with ${\bf x}_1={\bf o}$. Now, let us observe the local structures of $P_{2m}+X$ at the vertices ${\bf v}\in V+X$.

Let $X^{\bf v}$ denote the subset of $X$ consisting of all points ${\bf x}_i$ such that
$${\bf v}\in \partial (P_{2m})+{\bf x}_i.$$
Since $P_{2m}+X$ is a multiple tiling, the set $X^{\bf v}$ can be divided into disjoint subsets $X^{\bf v}_1$, $X^{\bf v}_2$, $\ldots ,$ $X^{\bf v}_t$ such that the translates in $P_{2m}+X^{\bf v}_j$ can be re-enumerated as $P_{2m}+{\bf x}^j_1$, $P_{2m}+{\bf x}^j_2$, $\ldots $, $P_{2m}+{\bf x}^j_{s_j}$ satisfying the following conditions:

\medskip
\noindent
{\bf 1.} {\it ${\bf v}\in \partial (P_{2m})+{\bf x}^j_i$ holds for all $i=1, 2, \ldots, s_j.$}

\smallskip\noindent
{\bf 2.} {\it Let $\angle^j_i$ denote the inner angle of $P_{2m}+{\bf x}^j_i$ at ${\bf v}$ with two half-line edges $L^j_{i,1}$ and $L^j_{i,2}$, where $L^j_{i,1}$, ${\bf x}^j_i-{\bf v}$ and $L^j_{i,2}$ are in clock order. Then, the inner angles join properly as
$$L^j_{i,2}=L^j_{i+1,1}$$
holds for all $i=1,$ $2,$ $\ldots ,$ $s_j$, where $L^j_{s_j+1,1}=L^j_{1,1}$.}

\medskip
For convenience, we call such a sequence $P_{2m}+{\bf x}^j_1$, $P_{2m}+{\bf x}^j_2$, $\ldots $, $P_{2m}+{\bf x}^j_{s_j}$ an {\it adjacent wheel} at ${\bf v}$. It is easy to see that
$$\sum_{i=1}^{s_j}\angle^j_i =2w_j\cdot \pi$$
hold for positive integers $w_j$. Then we define
$$\phi ({\bf v})=\sum_{j=1}^tw_j= {1\over {2\pi }}\sum_{j=1}^t\sum_{i=1}^{s_j}\angle^j_i$$
and
$$\varphi ({\bf v})=\sharp \left\{ {\bf x}_i:\ {\bf x}_i\in X,\ {\bf v}\in {\rm int}(P_{2m})+{\bf x}_i\right\}.$$

\medskip
Clearly, if $P_{2m}+X$ is a $\tau (P_{2m})$-fold translative tiling of $\mathbb{E}^2$, then
$$\tau (P_{2m})= \varphi ({\bf v})+\phi ({\bf v})\eqno (1)$$
holds for all ${\bf v}\in V+X$.

\medskip
First, let us introduce some basic results which will be useful in this paper.

\medskip\noindent
{\bf Lemma 1 (Bolle \cite{boll}).} {\it A convex polygon is a $k$-fold lattice tile for a lattice $\Lambda$ and some positive integer $k$ if and only if the following conditions are satisfied:

\noindent
{\bf 1.} It is centrally symmetric.

\noindent
{\bf 2.} When it is centered at the origin, in the relative interior of each edge $G$ there is a point of ${1\over 2}\Lambda $.

\noindent
{\bf 3.} If the midpoint of $G$ is not in ${1\over 2}\Lambda $ then $G$ is a lattice vector of $\Lambda $.}

\medskip\noindent
{\bf Lemma 2 (Zong \cite{zong}).} {\it If $m$ is even and $P_{2m}+\Lambda $ is a multiple lattice tiling, then $P_{2m}$ has an edge $G$ which is a lattice vector of $\Lambda $.}

\medskip\noindent
{\bf Lemma 3 (Zong \cite{zong}).} {\it Let ${\bf u}_i$ be the middle point of $G_i$. If $m$ is an odd positive integer, $P_{2m}+\Lambda $ is a $k$-fold lattice tiling of $\mathbb{E}^2$, and all ${\bf u}_i$ belong to ${1\over 2}\Lambda $, then we have
$$\sum_{i=1}^m(-1)^i{\bf u}_i ={\bf o},$$
where ${\bf o}=(0,0)$ is the origin of $\mathbb{E}^2$.}

\medskip\noindent
{\bf Lemma 4 (Yang and Zong \cite{yz2}).} {\it Assume that $P_{2m}$ is a centrally symmetric convex $2m$-gon centered at the origin and $P_{2m}+X$ is a $\tau (P_{2m})$-fold translative tiling of the plane, where $m\ge 4$. If ${\bf v}\in V+X$ is a vertex and $G\in \Gamma +X$ is an edge with ${\bf v}$ as one of its two ends, then there are at least $\lceil (m-3)/2\rceil $ different translates $P_{2m}+{\bf x}_i$ satisfying both
$${\bf v}\in \partial (P_{2m})+{\bf x}_i$$
and}
$$G\setminus \{ {\bf v}\}\subset {\rm int}(P_{2m})+{\bf x}_i.$$

\medskip\noindent
{\bf Lemma 5 (Yang and Zong \cite{yz2}).} {\it Let $P_{2m}$ be a centrally symmetric convex $2m$-gon, then}
$$\tau^*(P_{2m})\ge \tau (P_{2m})\ge \left\{\begin{array}{ll}
m-1,&\mbox{if $m$ is even,}\\
m-2,&\mbox{if $m$ is odd.}
\end{array}
\right.$$

\medskip\noindent
{\bf Lemma 6 (Yang and Zong \cite{yz2}).} {\it Assume that $P_{2m}$ is a centrally symmetric convex $2m$-gon centered at the origin, $P_{2m}+X$ is a translative multiple tiling of the plane, and $\mathbf{v}\in V+X$. Then we have
$$\phi ({\bf v})=\kappa\cdot {{m-1}\over 2}+\ell\cdot {1\over 2},$$
where $\kappa $ is a positive integer and $\ell$ is the number of the edges in $\Gamma+X$ which take ${\bf v}$ as an interior point.}

\vspace{0.6cm}
\noindent
{\Large\bf 3. Technical Lemmas}

\bigskip\noindent
{\bf Lemma 7.} {\it Let $P_{14}$ be a centrally symmetric convex tetradecagon, then}
$$\tau^* (P_{14})\ge \tau (P_{14})\ge 7.$$

\medskip
\noindent
{\bf Proof.} We take ${\bf v}\in V+X$ and assume that $P_{14}+{\bf x}_1$, $P_{14}+{\bf x}_2$, $\ldots $, $P_{14}+{\bf x}_s$ is an adjacent wheel at ${\bf v}$. First, it follows from Lemma 4  and Lemma 6 that
$$\varphi ({\bf v})\ge 2\eqno(2)$$
and
$$\phi ({\bf v})\ge 3.\eqno(3)$$

Now, we consider three cases.

\smallskip\noindent
{\bf Case 1.} {\it $\phi ({\bf v})\ge 5$ holds for a vertex ${\bf v}\in V+X$.} Then, by (1) and (2) we get
$$\tau (P_{14})= \varphi ({\bf v}) + \phi ({\bf v}) \ge 7.\eqno(4)$$

\smallskip\noindent
{\bf Case 2.} {\it $\phi ({\bf v})= 4$ holds for a vertex ${\bf v}\in V+X$.} Then, by Lemma 6 we get $\ell \not=0$. If ${\bf v}\in {\rm int}G$ holds for a suitable edge $G$, applying Lemma 4 to $G$ and its two ends we get
$$\varphi ({\bf v})\ge 4.$$
Then it follows by (1) that
$$\tau (P_{14})= \varphi ({\bf v}) + \phi ({\bf v}) \ge 8.\eqno(5)$$

\smallskip\noindent
{\bf Case 3.} {\it $\phi ({\bf v})= 3$ holds for every vertex ${\bf v}\in V+X$.} Then, the adjacent wheels at all ${\bf v}\in V$ are essentially unique, as shown by Figure 1. Let ${\bf v}_1$, ${\bf v}_2$, $\ldots $, ${\bf v}_{14}$ be the fourteen vertices of $P_{14}$. It follows that there are five point ${\bf y}_i\in X$ such that $P_{14}+{\bf x}_1$, $P_{14}+{\bf x}_7$, $P_{14}+{\bf y}_1$, $\ldots $, $P_{14}+{\bf y}_5$ is the adjacent wheel at ${\bf v}^*_1$. Then we have ${\bf v}_{10}+{\bf y}_2={\bf v}^*_1$, ${\bf v}_8+{\bf y}_4={\bf v}^*_1$ and
$${\bf v}\in {\rm int}(P_{14})+{\bf y}_i\quad i=2,\ 4.$$
By convexity, it can be easily deduced that
$${\bf v}^*_4\in {\rm int}(P_{14})+{\bf y}_i,\quad i=2,\ 4.$$

On the other hand, the adjacent wheel at ${\bf v}^*_4$ has two different translates taking ${\bf v}$ as an interior point as well. Thus, we have
$$\varphi ({\bf v})\ge 4$$
and
$$\tau (P_{14})= \varphi ({\bf v}) + \phi ({\bf v}) \ge 7.\eqno(6)$$

\begin{figure}[!ht]
\centering
\includegraphics[scale=0.43]{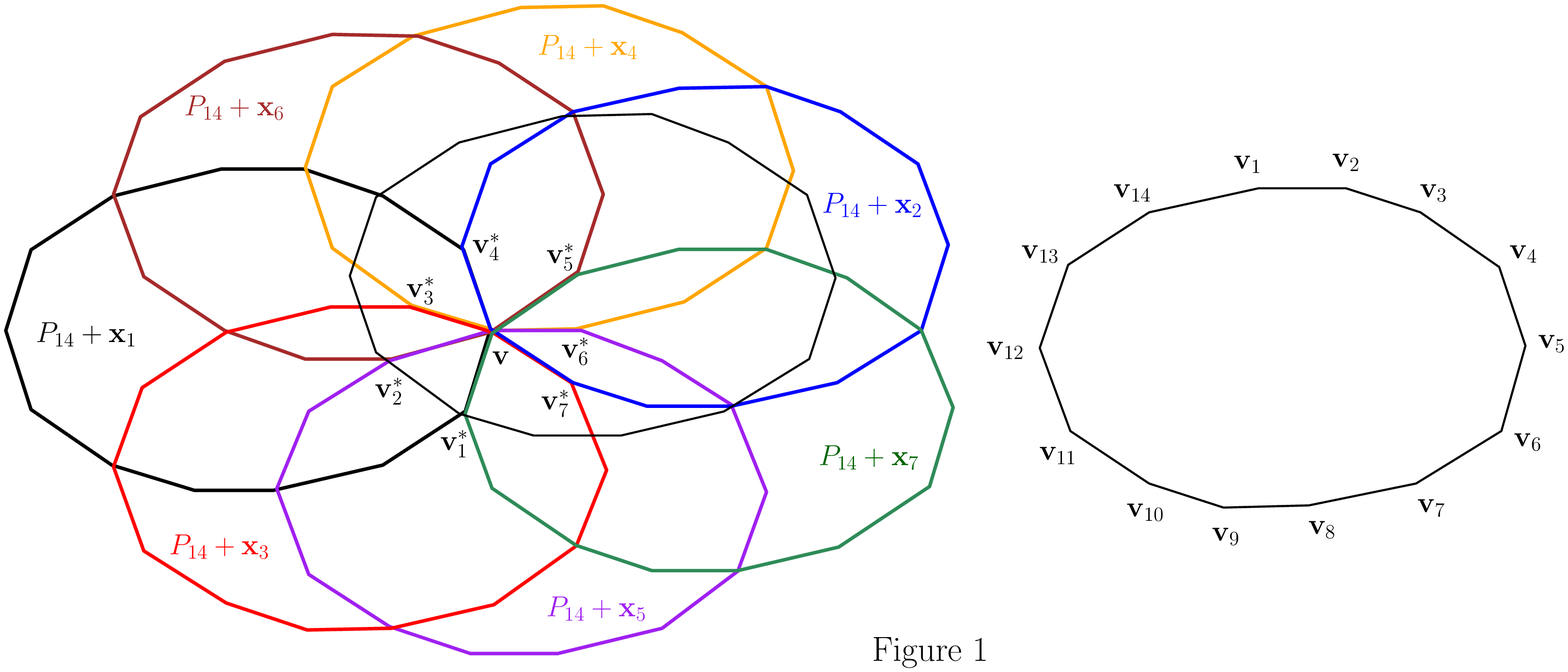}
\end{figure}

The lemma is proved.\hfill{$\Box$}

\medskip
\noindent
{\bf Lemma 8.} {\it Let $P_{12}$ be a centrally symmetric convex dodecagon, then}
$$\tau^* (P_{12})\ge 7.$$

\medskip
\noindent
{\bf Proof.} Since $\tau^*(P_{2m})$ is invariant under linear transformations on $P_{2m}$, we assume that $\Lambda =\mathbb{Z}^2$ and $P_{12}+\Lambda $ is a $\tau^*(P_{12})$-fold lattice tiling. Let ${\bf u}_i$ denote the middle point of $G_i$ and write ${\bf v}_i=(x_i,y_i)$ and ${\bf u}_i=(x'_i, y'_i)$. By Lemma 2 and a uni-modular transformation, we may assume that ${\bf v}_2-{\bf v}_1=(k, 0)$ and $y'_1>0$, where $k$ is a positive integer. By reduction (as shown by Figure 5), we many assume further that ${\bf v}_2-{\bf v}_1=(1,0)$. For convenience, let $P$ denote the parallelogram with vertices ${\bf v}_1$, ${\bf v}_2$, ${\bf v}_7=-{\bf v}_1$ and ${\bf v}_8=-{\bf v}_2$.

By Lemma 1 it follows that all $y_2-y_3$, $y_3-y_4$, $y_4-y_5$, $y_5-y_6$ and $y_6-y_7$ are positive integers. Thus, we have
$$y_1=y'_1=y_2\ge {5\over 2}.$$
If $y_1=y'_1=y_2\ge 3,$ then we have
$$\tau^*(P_{12})={\rm vol}(P_{12}) > {\rm vol}(P)\ge 6.\eqno(7)$$
If $y_1=y'_1=y_2= {5\over 2},$ then all ${\bf u}_i$ belong to ${1\over 2}\Lambda $. Let $T_i$ denote the triangle with vertices ${\bf u}_i$, ${\bf u}_{i+1}$ and ${\bf u}_6$, where $i=2, 3$ and $4$.
Clearly, all $y'_i-y'_6$ are positive integers. Thus, we have
$${\rm vol}(T_i)={1\over 2}\left|
\begin{array}{cc}
x'_i-x'_6& y'_i-y'_6\\
x'_{i+1}-x'_6& y'_{i+1}-y'_6
\end{array}\right|\ge {1\over 4}$$
and
$$\tau^*(P_{12})={\rm vol}(P_{12})>{\rm vol}(P)+2({\rm vol}(T_2)+{\rm vol}(T_3)+{\rm vol}(T_4) )\ge 5+6\cdot {1\over 4}>6.\eqno(8)$$

The Lemma is proved. \hfill{$\Box$}

\medskip\noindent
{\bf Lemma 9.} {\it Let $P_{10}$ be a centrally symmetric convex decagon, then
$$\tau^*(P_{10})=6$$
holds if and only if, under a suitable affine linear transformation, it takes ${\bf u}_1=(-1, {1\over 2})$, ${\bf u}_2=({1\over 2}, 1)$, ${\bf u}_3=({3\over 2}, 1)$, ${\bf u}_4=(2, {1\over 2})$, ${\bf u}_5=(2,0)$, ${\bf u}_6=-{\bf u}_1$, ${\bf u}_7=-{\bf u}_2$, ${\bf u}_8=-{\bf u}_3$, ${\bf u}_9=-{\bf u}_4$ and ${\bf u}_{10}=-{\bf u}_5$ as the middle points of their edges, or
takes ${\bf u}_1=({1\over 2}, -1)$, ${\bf u}_2=({3\over 2}, -{1\over 2})$, ${\bf u}_3=(2, 0)$, ${\bf u}_4=({3\over 2}, {1\over 2})$, ${\bf u}_5=({1\over 2}, 1)$, ${\bf u}_6=-{\bf u}_1$, ${\bf u}_7=-{\bf u}_2$, ${\bf u}_8=-{\bf u}_3$, ${\bf u}_9=-{\bf u}_4$ and ${\bf u}_{10}=-{\bf u}_5$ as the middle points of their edges.}

\medskip\noindent
{\bf Proof.} Let ${\bf v}_1$, ${\bf v}_2$, $\ldots $, ${\bf v}_{10}$ be the ten vertices of $P_{10}$ enumerated clock-wise, let $G_i$ denote the edge of $P_{10}$ with ends ${\bf v}_i$ and ${\bf v}_{i+1}$, where ${\bf v}_{11}={\bf v}_1$, and let ${\bf u}_i$ denote the middle point of $G_i$. For convenience, we write ${\bf v}_i=(x_i, y_i)$ and ${\bf u}_i=(x'_i, y'_i).$

It is known that $\sigma (D)+\sigma (\Lambda )$ is a $k$-fold lattice tiling of $\mathbb{E}^2$ whenever $D+\Lambda $ is such a tiling and $\sigma $ is a non-singular linear transformation from $\mathbb{E}^2$ to $\mathbb{E}^2$. Therefore, without loss of generality, we assume that $\Lambda =\mathbb{Z}^2$ and $P_{10}+\Lambda $ is a six-fold lattice tiling of $\mathbb{E}^2$.

By Lemma 1 we know that
$${\rm int}(G_i)\cap \mbox{${1\over 2}$}\Lambda \not=\emptyset $$
holds for all the ten edges $G_i$ and, if ${\bf u}_i\not \in {1\over 2}\Lambda $, then $G_i$ is a lattice vector of $\Lambda $. Now, we consider two cases.

\medskip\noindent
{\bf Case 1.} {\it $G_1$ is a lattice vector of $\Lambda $.} Without loss of generality, by a uni-modular linear transformation, we assume that ${\bf v}_2-{\bf v}_1=(k, 0)$ and $y'_1>0$, where $k$ is a positive integer. In fact, by reduction (as shown by Figure 5), one may assume that $G_1$ is primitive as a lattice vector and therefore $k=1$. Then, it can be deduced that
$$y_1=y'_1=y_2\in \mbox{${1\over 2}$}\mathbb{Z}$$
and all $y_i-y_{i+1}$ are integers. In particular, when $i=2$, $3,$ $4$ and $5$, they are positive integers. Thus, one can deduce that
$$y'_1= 2\ {\rm or}\ {5\over 2}.$$

\noindent
{\bf Case 1.1.} $y'_1=2$. Then we must have
$$y_2-y_3=y_3-y_4=y_4-y_5=y_5-y_6=1.$$
By the second term of Lemma 1, one can deduce that
$${\bf u}_i\in \mbox{${1\over 2}$} \Lambda,\quad i=2,\ 3,\ 4,\ 5. $$
Since ${\bf v}_2=(1,0)+{\bf v}_1$ and
$${\bf v}_{i+1}=2{\bf u}_i-{\bf v}_i, \quad i=2,\ 3,\ 4,\ 5,$$
it can be deduced that
$$-{\bf v}_1={\bf v}_6=2({\bf u}_5-{\bf u}_4+{\bf u}_3-{\bf u}_2)+(1,0)+{\bf v}_1$$
and therefore
$${\bf v}_i\in \mbox{${1\over 2}$} \Lambda, \quad i=1,\ 2,\ \ldots,\ 10. $$
Then all $G_i$ are lattice vectors.

\begin{figure}[!ht]
\centering
\includegraphics[scale=0.5]{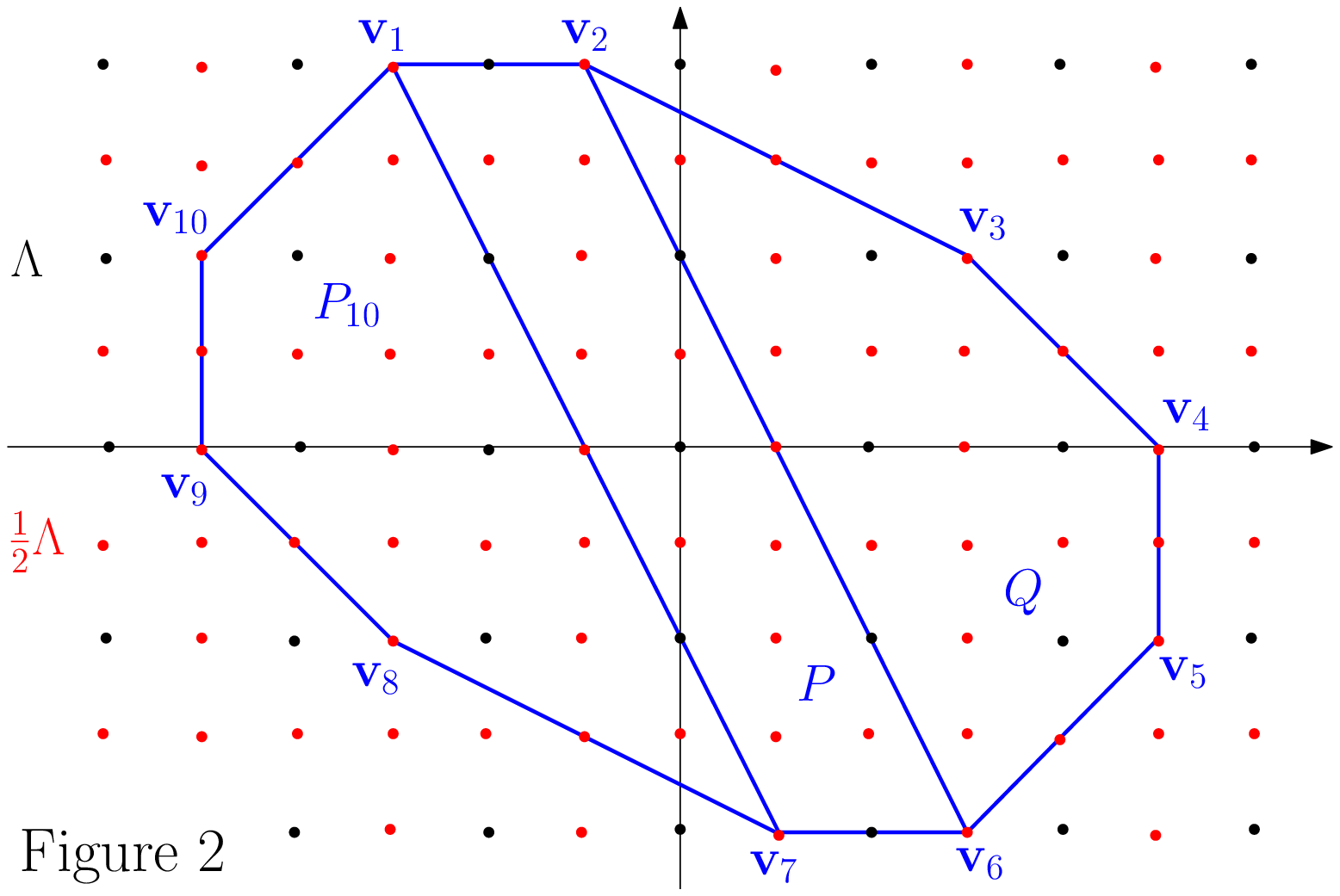}
\end{figure}

Let $P$ denote the parallelogram with vertices ${\bf v}_1$, ${\bf v}_2$, ${\bf v}_6$ and ${\bf v}_7$, and let $Q$ denote the pentagon with vertices ${\bf v}_2$, ${\bf v}_3$, ${\bf v}_4$, ${\bf v}_5$ and ${\bf v}_6$, as shown by Figure 2. Applying Pick's theorem to $Q$, we get
$${\rm vol}(Q)\ge \left({5\over 2}-1\right)$$
and therefore
$$\tau^*(P_{10})={\rm vol}(P_{10})={\rm vol}(P)+2\cdot {\rm vol}(Q)\ge 4+2\cdot \left( {5\over 2}-1\right)= 7.\eqno(9)$$

\noindent
{\bf Case 1.2.} $y'_1= {5\over 2}$. Then all $y_i-y_{i+1}$ are positive integers for $2\le i\le 5$. If
$${\bf u}_i\in \mbox{${1\over 2}$} \Lambda $$
hold for all $i=2,$ $3,$ $4$ and $5$, similar to the previous case one can deduce
$$\tau^*(P_{10})={\rm vol}(P_{10})\ge 7.\eqno(10)$$

If ${\bf u}_i\not\in {1\over 2}\Lambda$ holds for one of these indices, then we have $y_i-y_{i+1}=2$. By a uni-modular transformation, we may assume that $-{7\over 4}\le x_1<{3\over 4}.$ Then we have ${\bf v}_2-{\bf v}_6=(x, 5)$, where $-{5\over 2}\le x< {5\over 2}$. If ${\bf v}_i-{\bf v}_{i+1}=(k,2)$ with $|k|\ge 2$, let $Q$ denote the pentagon with vertices ${\bf v}_2$, ${\bf v}_3$, ${\bf v}_4$, ${\bf v}_5$ and ${\bf v}_6$, then we have
$${\rm vol}(Q)> {1\over 2}\left|
\begin{array}{cc}
x & 5\\
k& 2
\end{array}\right| = {1\over 2}\left| 2x-5k\right|\ge {5\over 2}$$
and thus
$$\tau^*(P_{10})={\rm vol}(P)+2\cdot {\rm vol}(Q)\ge 10.\eqno(11)$$
If ${\bf v}_i-{\bf v}_{i+1}=(k,2)$ with $k=\pm 1$, then we have $x_1\in {1\over 4}\mathbb{Z}$ and therefore $x\in {1\over 2}\mathbb{Z}$ and $-{5\over 2}\le x\le 2$. By considering two subcases with respect to $x_1=-{7\over 4}$ and $x_1\not=-{7\over 4}$, we can get
$${\rm vol} (Q)>{1\over 2}$$
and
$$\tau^*(P_{10})={\rm vol}(P)+2\cdot {\rm vol}(Q)>6.\eqno(12)$$

\medskip\noindent
{\bf Case 2.} {\it All the middle points ${\bf u}_i$ belong to ${1\over 2}\Lambda $.} Since $P_{10}+\Lambda $ is a six-fold lattice tiling of $\mathbb{E}^2$, one can deduce that
$${\rm vol}(2P_{10})=24$$
and all ${\bf u}'_i=2{\bf u}_i$ belong to $\Lambda $. For convenience, we define $Q_{10}$ to be the centrally symmetric lattice decagon with vertices ${\bf u}'_1$, ${\bf u}'_2$, $\ldots $, ${\bf u}'_{10}$ and write ${\bf u}'_i=(x'_i, y'_i)$. Since $Q_{10}$ is a centrally symmetric lattice polygon, its area must be a positive integer. Thus, we have
$${\rm vol}(Q_{10})\le 23.\eqno(13)$$

\medskip
Now, we explore $Q_{10}$ in detail by considering the following subcases.

\medskip\noindent
{\bf Case 2.1.} {\it ${\bf u}'_1$ is primitive in $\Lambda $}. Without loss of generality, guaranteed by uni-modular linear transformations, we take ${\bf u}'_1=(0, 1)$.  Then, Lemma 3 implies
$$\left\{
\begin{array}{ll}
x'_4-x'_5&\hspace{-0.3cm}=x'_3-x'_2, \\
y'_4-y'_5&\hspace{-0.3cm}=y'_3-y'_2+1.
\end{array}\right. \eqno(14) $$

If $x'_2\ge x'_3$ or $x'_3=x'_4$, one can easily deduce contradiction with convexity from (14). For example, if $x'_3=x'_4>x'_2$, then it can be deduced by (14) that
$${\bf u}'_2-{\bf u}'_5={\bf u}'_{10}-{\bf u}'_7=k{\bf u}'_1$$
with $k\ge 2$, which contradicts the assumption that $Q_{10}$ is a centrally symmetric convex decagon. Therefore, we may assume that
$$x'_3> x'_i\eqno(15)$$
for all $i\not= 3.$

Let $T'$ denote the lattice triangle with vertices ${\bf u}'_1$, ${\bf u}'_2$ and ${\bf u}'_3$, let $Q$ denote the lattice quadrilateral with vertices ${\bf u}'_3$, ${\bf u}'_4$, ${\bf u}'_5$ and ${\bf u}'_6$, and let $T$ denote the lattice triangle with vertices ${\bf u}'_1$, ${\bf u}'_3$ and ${\bf u}'_6$ (as shown by Figure 3). It follows from (13) and Pick's theorem that
$${\rm vol}(T)\le {1\over 2}\Bigl(23 -2 \bigl({\rm vol}(T')+{\rm vol}(Q)\bigr)\Bigr)\le 10$$
and therefore
$$x'_3\le 10.\eqno(16)$$

\begin{figure}[!ht]
\centering
\includegraphics[scale=0.6]{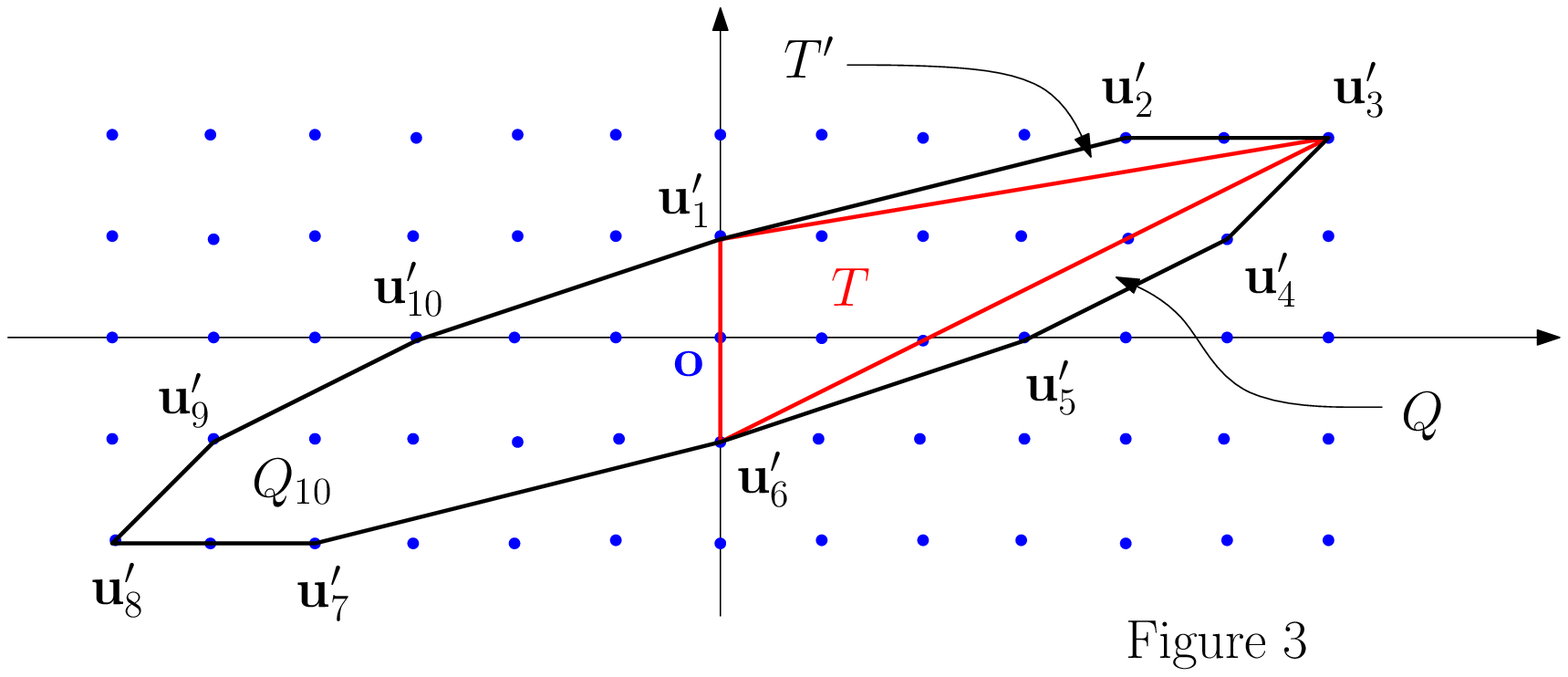}
\end{figure}

Let $\alpha $ denote the slope of $G_1$, that is
$$\alpha ={{y_2-y_1}\over {x_2-x_1}}.$$
By a uni-modular linear transformation such as
$$\left\{
\begin{array}{ll}
x'=x,&\\
y'=y+kx,&
\end{array}
\right.$$
where $k$ is a suitable integer, we may assume that
$$0\le \alpha <1.\eqno (17)$$
Let $L_i$ denote the straight line containing $G_i$, it is obvious that $P_{10}$ is in the strip bounded by $L_1$ and $L_6$. Furthermore, we define
five slopes
$$\beta_i={{y'_{i+1}-y'_i}\over {x'_{i+1}-x'_i}},\quad i=1,\ 2,\ \ldots,\ 5.$$

By convexity it can be shown that there is no six-fold lattice decagon tile with $\alpha =0$ in our setting. When $\alpha >0$, by (14) and convexity it follows that $y'_4-y'_5\ge 1$ and therefore
$$y'_3-y'_2\ge 0.\eqno (18)$$

\begin{figure}[!ht]
\centering
\includegraphics[scale=0.6]{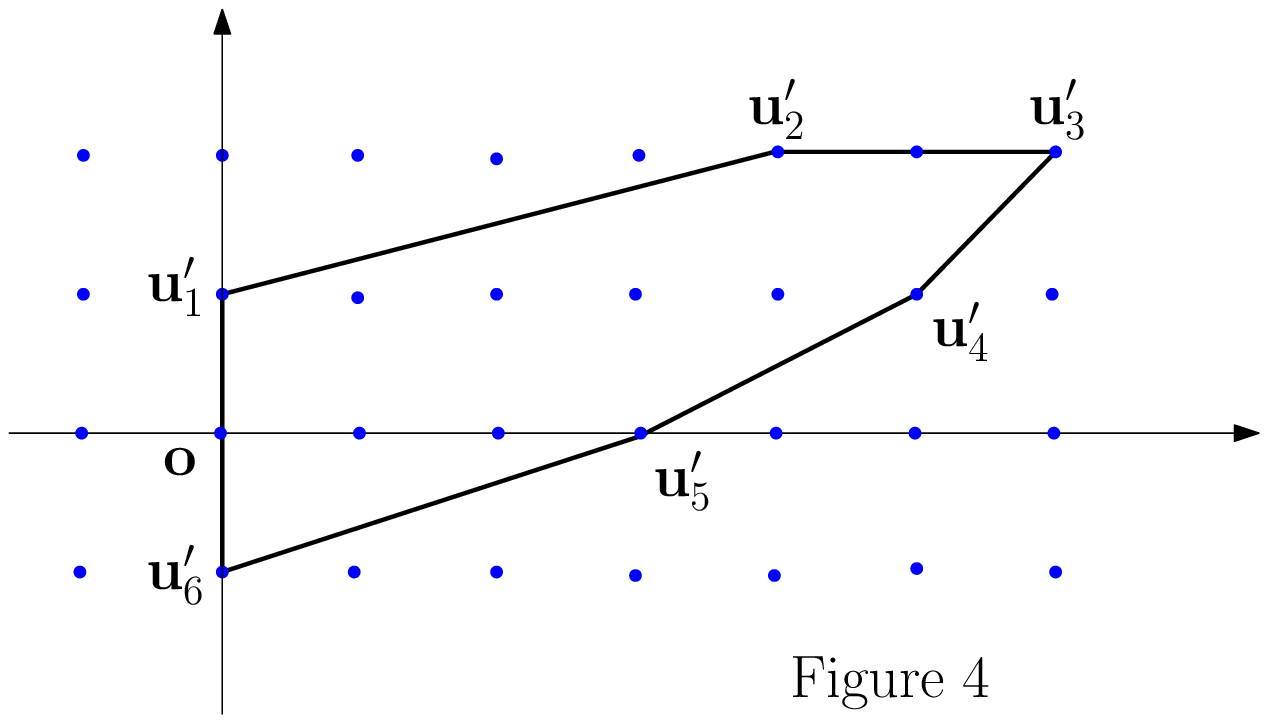}
\end{figure}

As shown by Figure 4, we assume that $${\bf u}'_3-{\bf u}'_4=(p_1, q_1)$$ and $${\bf u}'_5-{\bf u}'_6=(p_2, q_2),$$
where all $p_i$ and $q_i$ are positive integers. Then, by (16) we have
$$x'_3-x'_2=x'_3-(x'_2-x'_1)=x'_3-(p_1+p_2)\le 8.$$
Now, we consider in subcases with respect to the different orientations of ${\bf u}'_3-{\bf u}'_2.$

\medskip\noindent
{\bf Case 2.1.1.} {\it $y'_3-y'_2=0$ and $x'_3-x'_2=1$}. By (14) and convexity we have $x'_4-x'_5=1$, $y'_4-y'_5=1$, $\beta_4=1$,
$$\beta_3={{q_1}\over {p_1}}>1$$
and
$$\beta_5={{q_2}\over {p_2}}<1.$$
Then, one can deduce that
$$\beta_1={{q_1+q_2-1}\over {p_1+p_2}}>{{q_2}\over {p_2}}=\beta_5,$$
which contradicts the convexity of $Q_{10}$.

\medskip\noindent
{\bf Case 2.1.2.} {\it $y'_3-y'_2=0$ and $x'_3-x'_2=2$}. By (14) and convexity we have $x'_4-x'_5=2$, $y'_4-y'_5=1$, $\beta_4={1\over 2}$,
$$\beta_3={{q_1}\over {p_1}}>{1\over 2},\eqno(19)$$
$$\beta_5={{q_2}\over {p_2}}<{1\over 2}\eqno(20)$$
and
$$\beta_1={{q_1+q_2-1}\over {p_1+p_2}}<{{q_2}\over {p_2}}.\eqno(21)$$

By (16) and (20) one can deduce that
$$3\le p_2\le 7,\eqno(22)$$
$$1\le p_1\le 5\eqno(23)$$
and
$$1\le q_2\le 3.\eqno(24)$$
On the other hand, by (21), (23) and (24) we get
$$q_1<p_1\cdot {{q_2}\over {p_2}}+1<{1\over 2}\cdot p_1+1$$
and therefore
$$1\le q_1\le 3.\eqno(25)$$

Then, it can be verified that the only integer groups $(p_1, q_1, p_2,q_2)$ satisfying (16), (19), (20) and (21) are $(1,1,3,1),$ $(1,1,4,1)$, $(1,1,5,1)$, $(1,1,6,1),$ $(1,1,7,1)$, $(1,1,5,2)$, $(1,1,6,2)$, $(1,1,7,2)$ and $(1,1,7,3)$. By checking the areas of their corresponding decagons, keeping the subcase conditions in mind, there are only two $Q_{10}$ satisfying (13). Namely, the one with vertices ${\bf u}'_1=(0,1)$, ${\bf u}'_2=(4,2)$, ${\bf u}'_3=(6,2)$, ${\bf u}'_4=(5,1)$, ${\bf u}'_5=(3,0)$, ${\bf u}'_6=-{\bf u}'_1$, ${\bf u}'_7=-{\bf u}'_2$, ${\bf u}'_8=-{\bf u}'_3$, ${\bf u}'_9=-{\bf u}'_4$ and ${\bf u}'_{10}=-{\bf u}'_5$, which produces the five-fold lattice tiles, and the one with vertices ${\bf u}'_1=(0,1)$, ${\bf u}'_2=(5,2)$, ${\bf u}'_3=(7,2)$, ${\bf u}'_4=(6,1)$, ${\bf u}'_5=(4,0)$, ${\bf u}'_6=-{\bf u}'_1$, ${\bf u}'_7=-{\bf u}'_2$, ${\bf u}'_8=-{\bf u}'_3$, ${\bf u}'_9=-{\bf u}'_4$ and ${\bf u}'_{10}=-{\bf u}'_5$, which indeed produces six-fold lattice tiles. Clearly, the second decagon is equivalent to the first one stated in the lemma under the linear transformation
$$\left\{
\begin{array}{ll}
x'={1\over 2}(x-2y),&\\
\vspace{-0.3cm}
&\\
y'={1\over 2}y.&
\end{array}
\right.$$

\medskip\noindent
{\bf Case 2.1.3.} {\it $y'_3-y'_2=0$ and $x'_3-x'_2=3$}. By (14) and convexity we have $x'_4-x'_5=3$, $y'_4-y'_5=1$, $\beta_4={1\over 3}$,
$$\beta_3={{q_1}\over {p_1}}>{1\over 3},\eqno(26)$$
$$\beta_5={{q_2}\over {p_2}}<{1\over 3}\eqno(27)$$
and
$$\beta_1={{q_1+q_2-1}\over {p_1+p_2}}<{{q_2}\over {p_2}}.\eqno(28)$$

Restricted by (16), similar to the previous case, it can be deduced that the only integer solutions $(p_1,q_1, p_2,q_2)$  for (26), (27) and (28) are $(1,1,4,1)$, $(1,1,5,1)$,
$(1,1,6,1)$, $(2,1,4,1)$, $(2,1,5,1)$ and $(2,1,6,1)$. Then one can deduce
$${\rm vol}(Q_{10})\ge 25 \eqno(29)$$
for all these cases, which contradicts (13).

\medskip\noindent
{\bf Case 2.1.4.} {\it $y'_3-y'_2=0$ and $x'_3-x'_2 = 4$}. Then, one can easily deduce that $\beta_4={1\over 4}$,
$$\beta_3={{q_1}\over {p_1}}>{1\over 4},$$
$$\beta_5={{q_2}\over {p_2}}<{1\over 4} $$
and
$$\beta_1={{q_1+q_2-1}\over {p_1+p_2}}<{{q_2}\over {p_2}}.$$
Restricted by (16), similar to Case 2.1.3 one can deduce that the only integer solutions $(p_1,q_1, p_2,q_2)$  for these inequalities are $(1,1,5,1)$, $(2,1,5,1)$,
and $(3,1,5,1)$. Then we have
$${\rm vol}(Q_{10})\ge 25\eqno(30)$$
for all these cases, which contradicts (13).

\medskip\noindent
{\bf Case 2.1.5.} {\it $y'_3-y'_2=0$ and $x'_3-x'_2 \ge 5$}. Then, one can easily deduce that $\beta_4\le {1\over 5}$, $p_2\ge 6$ and
$$x'_3\ge 5+6>10,$$
which contradicts the restriction of $(16)$.

\medskip\noindent
{\bf Case 2.1.6.} {\it $y'_3-y'_2=1$ and $x'_3-x'_2=1$}. Then, by convexity we get
$$\alpha >\beta_1>\beta_2=1,\eqno(31)$$
which contradicts the assumption of (17).

\medskip\noindent
{\bf Case 2.1.7.} {\it $y'_3-y'_2=1$ and $x'_3-x'_2=2$}. By (14) and convexity we get $x'_4-x'_5=2$, $y'_4-y'_5=2,$ $\beta_4=1$,
$$\beta_3={{q_1}\over {p_1}}>1$$
and
$$\beta_5={{q_2}\over {p_2}}<1.$$
Then, it can be deduced that
$$\beta_1={{q_1+q_2-1}\over {p_1+p_2}}> {{q_2}\over {p_2}}=\beta_5,\eqno(32)$$
which contradicts the convexity assumption of $Q_{10}$.

\medskip\noindent
{\bf Case 2.1.8.} {\it $y'_3-y'_2=1$ and $x'_3-x'_2=3$}. Then we have $x'_4-x'_5=3$, $y'_4-y'_5=2,$ $\beta_2={1\over 3}$ and $\beta_4={2\over 3}$.

On one hand, by (16) it follows that $p_2\le 6$. On the other hand, by $\beta_2<\beta_1<\beta_5<\beta_4$ it follows that
$${1\over 3}<{{q_2}\over {p_2}}<{2\over 3}.$$
Thus, the integer pair $(p_2, q_2)$ has only five choices $(2,1)$, $(4,2)$, $(5,2)$, $(5,3)$ and $(6,3)$. Then, by checking
$${{q_1}\over {p_1}}>{2\over 3},$$
$${1\over 3}<{{q_1+q_2-1}\over {p_1+p_2}}<{{q_2}\over {p_2}}$$
and
$$p_1+p_2\le 7,$$
it can be deduced that the only candidates for $(p_1,q_1, p_2,q_2)$ are $(1,1,4,2)$, $(1,1,5,3)$, $(2,2,5,3)$ and $(1,1,6,3)$. In fact, the only candidate satisfying (13) is the one with vertices ${\bf u}'_1=(0,1)$, ${\bf u}'_2=(5,3)$, ${\bf u}'_3=(8,4)$, ${\bf u}'_4=(7,3)$, ${\bf u}'_5=(4,1)$, ${\bf u}'_6=-{\bf u}'_1$, ${\bf u}'_7=-{\bf u}'_2$, ${\bf u}'_8=-{\bf u}'_3$, ${\bf u}'_9=-{\bf u}'_4$ and ${\bf u}'_{10}=-{\bf u}'_5$, satisfying
$${\rm vol}(Q_{10})=22.\eqno(33)$$
This decagon indeed produces six-fold lattice tiles. Clearly, it is equivalent to the second one stated in the lemma under the linear transformation
$$\left\{
\begin{array}{ll}
x'={1\over 2}y,&\\
\vspace{-0.3cm}
&\\
y'={1\over 2}(x-2y).&
\end{array}
\right.$$

\medskip\noindent
{\bf Case 2.1.9.} {\it $y'_3-y'_2=1$ and $x'_3-x'_2=4$}. By (14), (16) and convexity it can be deduced that $p_2\le 5$, $\beta_4={1\over 2}$ and $\beta_5<\beta_4.$ Consequently, we have $\beta_5={1\over 3},$ ${1\over 4},$ ${1\over 5}$ or ${2\over 5}.$ Thus, by $\beta_2={1\over 4}$ and $\beta_2<\beta_1<\beta_5$ we get
$${1\over 4}<{{q_1+q_2-1}\over {p_1+p_2}}<{2\over 5}.\eqno (34)$$
By (16) we have $p_1+p_2\le 6$ and therefore (34) has only one solution $(p_1,q_1,p_2,q_2)=(1,1,5,2)$. However, for such decagon we have
$${\rm vol}(Q_{10})=29,\eqno(35)$$
which contradicts $(13)$.

\medskip\noindent
{\bf Case 2.1.10.} {\it $y'_3-y'_2=1$ and $x'_3-x'_2= 5$}. Then by (16) and convexity we have
$$p_1+p_2\le 5$$
and
$${1\over 5}<{{q_1+q_2-1}\over {p_1+p_2}}<{2\over 5}.$$
In fact, these inequalities have no positive integer solution.

\medskip\noindent
{\bf Case 2.1.11.} {\it $y'_3-y'_2=1$ and $x'_3-x'_2\ge 6$}. It follows by (16) that $p_2\le 3$. Then we get both $\beta_4\le {1\over 3}$ and $\beta_5\ge {1\over 3}$, which contradicts the convexity of $Q_{10}$.

\medskip\noindent
{\bf Case 2.1.12.} {\it $y'_3-y'_2=2$ and $x'_3-x'_2=3$}. Then by (14) and convexity we get $\beta_4=1$ and $\beta_1<\beta_5$. However
the two inequalities
$${{q_1}\over {p_1}}>\beta_4=1$$
and
$${{q_1+q_2-1}\over {p_1+p_2}}<{{q_2}\over {p_2}}$$
have no integer solution.

\medskip\noindent
{\bf Case 2.1.13.} {\it $y'_3-y'_2=2$ and $x'_3-x'_2=4$}. Then by (14) and convexity we get $\beta_2={1\over 2}$, $\beta_4={3\over 4}$, $\beta_2<\beta_5<\beta_4$ and therefore
$${1\over 2}<{{q_2}\over {p_2}}<{3\over 4}.\eqno(36)$$
Clearly, by (16) we have $p_2\le 5$ and therefore (36) has two groups of integer solutions $(p_2,q_2)=(3,2)$ or $(5,3)$. Then, the two inequalities
$p_1+p_2\le 6$ and
$${1\over 2}<{{q_1+q_2-1}\over {p_1+p_2}}<{{q_2}\over {p_2}}$$
have one group of integer solution $(p_1,q_1,p_2,q_2)=(2,2,3,2)$. Unfortunately, then we have
$${\rm vol}(Q_{10})=25,\eqno(37)$$
which contradicts (13).

\medskip\noindent
{\bf Case 2.1.14.} {\it $y'_3-y'_2=2$ and $x'_3-x'_2=5$}. Then by (14) and convexity we get $\beta_2={2\over 5}$, $\beta_4={3\over 5}$, $\beta_2<\beta_5<\beta_4$ and therefore
$${2\over 5}<{{q_2}\over {p_2}}<{3\over 5}.\eqno(38)$$
Clearly, by (16) we have $p_2\le 4$ and therefore (38) has two groups of integer solutions $(p_2,q_2)=(2,1)$ or $(4,2)$. Then, one can deduce that
$p_1+p_2\le 5$ and
$${2\over 5}<{{q_1+q_2-1}\over {p_1+p_2}}<{{q_2}\over {p_2}}$$
have no integer solution.

\medskip\noindent
{\bf Case 2.1.15.} {\it $y'_3-y'_2=2$ and $x'_3-x'_2=6$}. Then by (14) and convexity we get $\beta_5<\beta_4={1\over 2}$ and therefore $\beta_5={1\over 3}$, which contradicts the fact
$$\beta_5>\beta_1>\beta_2={1\over 3}.\eqno(39)$$

\medskip\noindent
{\bf Case 2.1.16.} {\it $y'_3-y'_2=2$ and $x'_3-x'_2\ge 7$}. Then by (14) and convexity we get $\beta_4\le {3\over 7}$ and $\beta_5\ge {1\over 2}$, which contradicts the convexity of $Q_{10}$.

\medskip\noindent
{\bf Case 2.1.17.} {\it $y'_3-y'_2=3$ and $x'_3-x'_2=4$}. Then by (14) and convexity we have $p_2\le 5$, $\beta_2={3\over 4}$ and $\beta_4=1$. Then we have
$$\beta_3={{q_1}\over {p_1}}>1$$
and therefore
$$\beta_1={{q_1+q_2-1}\over {p_1+p_2}}>{{q_2}\over {p_2}}=\beta_5,\eqno(40)$$
which contradicts the convexity of $Q_{10}$.

\medskip\noindent
{\bf Case 2.1.18.} {\it $y'_3-y'_2=3$ and $x'_3-x'_2=5$}. Then by (14) and convexity we have $p_2\le 4$, $\beta_2={3\over 5},$ $\beta_4={4\over 5}$ and $\beta_2<\beta_5<\beta_4$. The inequalities $p_2\le 4$ and
$${3\over 5}<{{q_2}\over {p_2}}<{4\over 5}$$
have two solutions $(p_2,q_2)=(3,2)$ or $(4,3)$. Then
$${3\over 5}<{{q_1+q_2-1}\over {p_1+p_2}}<{{q_2}\over {p_2}}\eqno(41)$$
has no solution satisfying $p_1+p_2\le 5$.

\medskip\noindent
{\bf Case 2.1.19.} {\it $y'_3-y'_2=3$ and $x'_3-x'_2=6$}. Then by (14) and convexity we get $p_2\le 3$, $\beta_2={1\over 2}$, $\beta_4={2\over 3}$ and $\beta_2< \beta_5<\beta_4$.  Then the inequalities $p_2\le 3$ and
$${1\over 2}<{{q_2}\over {p_2}}<{2\over 3}$$
have no solution.

\medskip\noindent
{\bf Case 2.1.20.} {\it $y'_3-y'_2=3$ and $x'_3-x'_2=7$}. Then by (14) and convexity we get $p_2\le 2$, $\beta_2={3\over 7}$, $\beta_4={4\over 7}$ and $\beta_2< \beta_5<\beta_4$.  Then the inequalities $p_2\le 2$ and
$${3\over 7}<{{q_2}\over {p_2}}<{4\over 7}$$
have one solution $(p_2,q_2)=(2,1)$. Then
$${3\over 7}<{{q_1}\over {p_1+2}}<{1\over 2}\eqno(42)$$
has no solution.

\medskip\noindent
{\bf Case 2.1.21.} {\it $y'_3-y'_2=3$ and $x'_3-x'_2\ge 8$}. Then by (14) and convexity we get $p_2=1$, $\beta_5\ge 1$ and $\beta_4\le {1\over 2}$, which contradicts the convexity of $Q_{10}$.

\medskip\noindent
{\bf Case 2.1.22.} {\it $y'_3-y'_2=4$ and $x'_3-x'_2=5$}. Then by (14) and convexity we have $p_2\le 4$, $\beta_2={4\over 5},$ $\beta_4=1$ and $\beta_2<\beta_5<\beta_4$. The inequalities $p_2\le 4$ and
$${4\over 5}<{{q_2}\over {p_2}}<1\eqno(43)$$
have no common integer solution.

\medskip\noindent
{\bf Case 2.1.23.} {\it $y'_3-y'_2=4$ and $x'_3-x'_2=6$}. Then by (14) and convexity we get $p_2\le 3$, $\beta_2={2\over 3},$ $\beta_4={5\over 6}$ and $\beta_2<\beta_5<\beta_4$. The inequalities $p_2\le 3$ and
$${2\over 3}<{{q_2}\over {p_2}}<{5\over 6}\eqno(44)$$
have no common integer solution.

\medskip\noindent
{\bf Case 2.1.24.} {\it $y'_3-y'_2=4$ and $x'_3-x'_2=7$}. Then by (14) and convexity we get $p_2\le 2$, $\beta_2={4\over 7},$ $\beta_4={5\over 7}$ and $\beta_2<\beta_5<\beta_4$. The inequalities $p_2\le 2$ and
$${4\over 7}<{{q_2}\over {p_2}}<{5\over 7}\eqno(45)$$
have no common integer solution.

\medskip\noindent
{\bf Case 2.1.25.} {\it $y'_3-y'_2=4$ and $x'_3-x'_2\ge 8$}. Then by (14) and convexity we get $p_2=1$, $\beta_5\ge 1$ and $\beta_4\le {5\over 8}$, which contradicts the convexity of $Q_{10}$.

\medskip\noindent
{\bf Case 2.1.26.} {\it $y'_3-y'_2=5$ and $x'_3-x'_2=6$}. Then by (14) and convexity we get $p_2\le 3$, $\beta_2={5\over 6},$ $\beta_4=1$ and $\beta_2<\beta_5<\beta_4$. The inequalities $p_2\le 3$ and
$${5\over 6}<{{q_2}\over {p_2}}<1\eqno(46)$$
have no common integer solution.

\medskip\noindent
{\bf Case 2.1.27.} {\it $y'_3-y'_2=5$ and $x'_3-x'_2=7$}. Then by (14) and convexity we get $p_2\le 2$, $\beta_2={5\over 7},$ $\beta_4={6\over 7}$ and $\beta_2<\beta_5<\beta_4$. The inequalities $p_2\le 2$ and
$${5\over 7}<{{q_2}\over {p_2}}<{6\over 7}\eqno(47)$$
have no common integer solution.

\medskip\noindent
{\bf Case 2.1.28.} {\it $y'_3-y'_2=5$ and $x'_3-x'_2\ge 8$}. Then by (14) and convexity we get $p_2=1$, $\beta_5\ge 1$ and $\beta_4\le {6\over 8}$, which contradicts the convexity of $Q_{10}$.

\medskip\noindent
{\bf Case 2.1.29.} {\it $y'_3-y'_2=6$ and $x'_3-x'_2=7$}. Then by (14) and convexity we get $p_2\le 2$, $\beta_2={6\over 7},$ $\beta_4=1$ and $\beta_2<\beta_5<\beta_4$. The inequalities $p_2\le 2$ and
$${6\over 7}<{{q_2}\over {p_2}}<1\eqno(48)$$
have no common integer solution.

\medskip\noindent
{\bf Case 2.1.30.} {\it $y'_3-y'_2=6$ and $x'_3-x'_2\ge 8$}. Then by (14) and convexity we get $p_2=1$, $\beta_5\ge 1$ and $\beta_4\le {7\over 8}$, which contradicts the convexity of $Q_{10}$.

\medskip\noindent
{\bf Case 2.1.31.} {\it $y'_3-y'_2=7$ and $x'_3-x'_2\ge 8$}. Then by (14) and convexity we get $p_2=1$, $\beta_5\ge 1$ and $\beta_4\le 1$, which contradicts the convexity of $Q_{10}$.

\medskip\noindent
{\bf Case 2.2.} {\it All ${\bf u}'_i$ are even multiplicative}. Then all ${\bf u}_i$ belong to $\Lambda $. It follows by Lemma 1 that ${1\over 2}P_{10}+\Lambda $ is a $k$-fold lattice tiling with
$$k={\rm vol}\left(\mbox{${1\over 2}$}P_{10}\right)={6\over 4}={3\over 2},\eqno(49)$$
which contradicts the fact that $k$ is a positive integer.

\medskip\noindent
{\bf Case 2.3.} {\it All ${\bf u}'_i$ are multiplicative, ${\bf u}'_1$ is odd multiplicative}. Without loss of generality, guaranteed by uni-modular linear transformations, we take ${\bf u}'_1=(0, 2q+1)$, where $q$ is a positive integer.

By Lemma 3 it follows that
$$x'_4-x'_5=x'_3-x'_2.$$
Therefore, by convexity and reflection we may assume that
$$x'_3\ge x'_i, \qquad i=1, 2, \ldots , 10.$$

Let $T'$ denote the lattice triangle with vertices ${\bf u}'_1$, ${\bf u}'_2$ and ${\bf u}'_3$, let $Q$ denote the lattice quadrilateral with vertices ${\bf u}'_3$, ${\bf u}'_4$, ${\bf u}'_5$ and ${\bf u}'_6$, and let $T$ denote the lattice triangle with vertices ${\bf u}'_1$, ${\bf u}'_3$ and ${\bf u}'_6$, as shown in Figure 3. It follows from (13) and Pick's theorem that
$${\rm vol}(T)\le {1\over 2}\Bigl( 23 -2 \bigl({\rm vol}(T')+{\rm vol}(Q)\bigr)\Bigr)\le 10\eqno(50)$$
and therefore
$$x'_3 ={{2\cdot {\rm vol}(T)}\over {2(2q+1)}}\le \left\lfloor {{10}\over 3}\right\rfloor = 3.\eqno(51)$$

It is assumed that all ${\bf u}'_i$ are multiplicative. Therefore by convexity we have
$$x'_2=x'_5=2$$
and
$$x'_3=x'_4=3.$$
Then, we have
$${\rm vol}(Q_{10})\ge 3\cdot (2(2q+1)+3)\ge 27,\eqno(52)$$
which contradicts (13).

\medskip
As a conclusion of all these cases, Lemma 9 is proved.\hfill{$\Box$}

\medskip\noindent
{\bf Lemma 10.} {\it Let $P_8$ be a centrally symmetric convex octagon, then
$$\tau^*(P_8)=6$$
if and only if (under suitable affine linear transformations) it with vertices ${\bf v}_1=(\alpha-1 , 2)$, ${\bf v}_2=(\alpha , -2)$, ${\bf v}_3=(1-\alpha , 0)$, ${\bf v}_4=(1+\alpha , -1)$, ${\bf v}_5=-{\bf v}_1$, ${\bf v}_6=-{\bf v}_2$, ${\bf v}_7=-{\bf v}_3$ and ${\bf v}_8=-{\bf v}_4$, where $0<\alpha <{1\over 6}$.}

\medskip\noindent
{\bf Proof.} Let $P_8$ be a centrally symmetric convex octagon centered at the origin, let ${\bf v}_1$, ${\bf v}_2$, $\ldots$, ${\bf v}_8$ be the eight vertices of $P_8$ enumerated in an anti-clock order, let $G_i$ denote the edge with ends ${\bf v}_i$ and ${\bf v}_{i+1}$, where ${\bf v}_9={\bf v}_1$, and let ${\bf u}_i$ denote the midpoint of $G_i$. For convenience, we write ${\bf v}_i=(x_i,y_i)$ and ${\bf u}_i=(x'_i,y'_i)$. Assume that $\Lambda = \mathbb{Z}^2$ and $P_8+\Lambda $ is a six-fold lattice tiling. Then, we have
$$\tau^*(P_8)={\rm vol}(P_8)=6.\eqno(53)$$

Based on Lemma 2, by a uni-modular transformation, we may assume that $G_1\cap {1\over 2}\Lambda\not=\emptyset$ and ${\bf v}_2-{\bf v}_1=(k, 0),$
where $k$ is a positive integer.  If $k>1$, we define $P_8'$ to be the octagon with vertices ${\bf v}'_1={\bf v}_1+({{k-1}\over 2},0),$ ${\bf v}'_2={\bf v}_2+({{1-k}\over 2},0),$ ${\bf v}'_3={\bf v}_3+({{1-k}\over 2},0),$ ${\bf v}'_4={\bf v}_4+({{1-k}\over 2},0),$ ${\bf v}'_5={\bf v}_5+({{1-k}\over 2},0),$ ${\bf v}'_6={\bf v}_6+({{k-1}\over 2},0),$ ${\bf v}'_7={\bf v}_7+({{k-1}\over 2},0)$ and ${\bf v}'_8={\bf v}_8+({{k-1}\over 2},0)$, as shown by Figure 5. By Lemma 1 it can be shown that
$P'_8+\Lambda $ is a multiple lattice tiling of $\mathbb{E}^2$ and therefore
$$\tau^*(P'_8)\le {\rm vol}(P'_8)\le {\rm vol}(P_8)-3=3,\eqno(54)$$
which contradicts the known fact that $\tau^* (P'_8)\ge 5$ (see \cite{yz1}). Thus, we have ${\bf v}_2-{\bf v}_1=(1,0).$

\begin{figure}[!ht]
\centering
\includegraphics[scale=0.52]{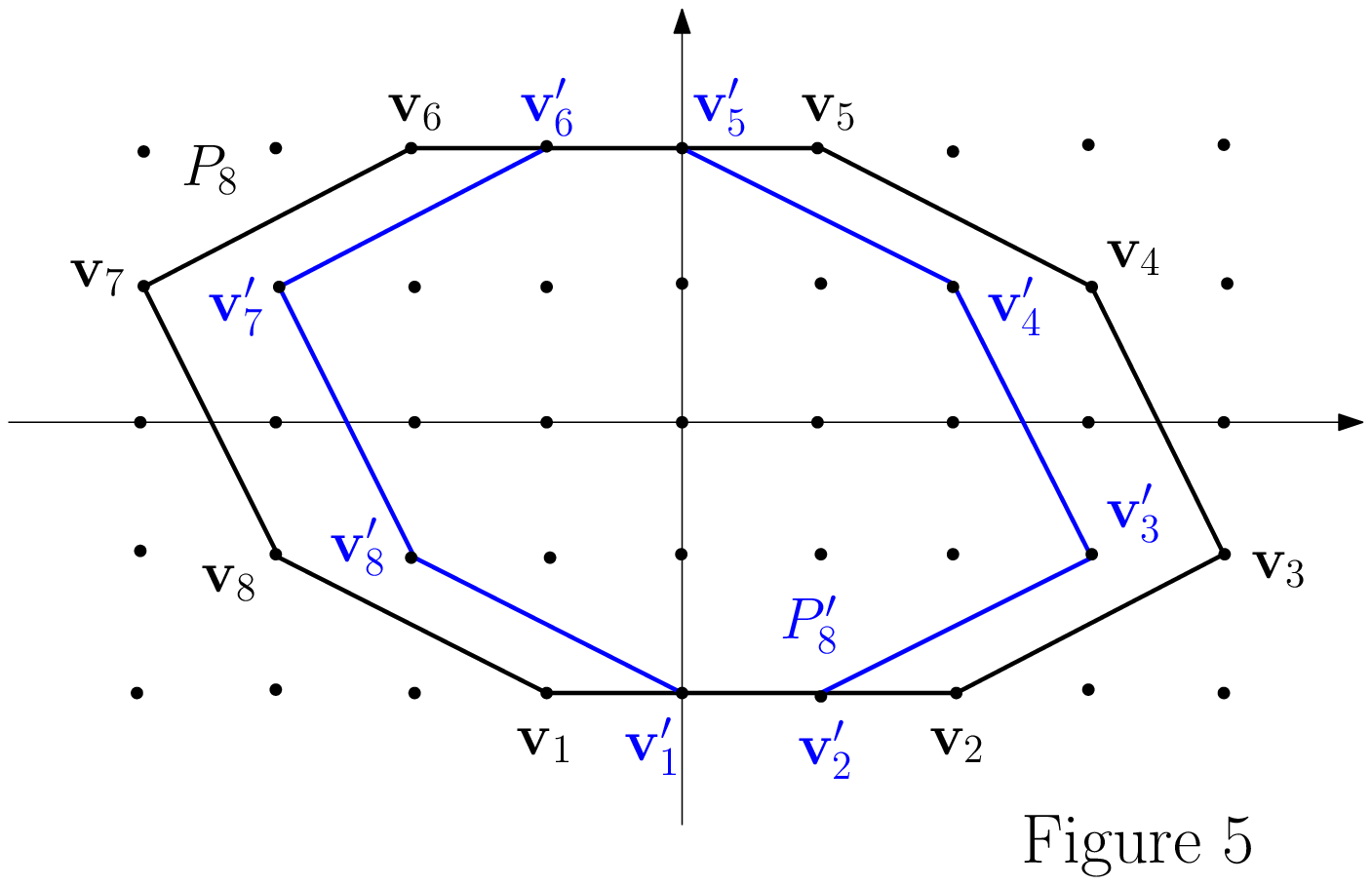}
\end{figure}

Apply Lemma 1 successively to $G_1$, $G_2$, $G_3$ and $G_4$, one can deduce that all $2y_2$, $y_3-y_2$, $y_4-y_3$ and $y_5-y_4$ are positive integers. Therefore, we have
$$y_2=y_1\le -{3\over 2}.\eqno(55)$$
On the other hand, if $y_2=y_1\le -3$ and let $P$ denote the parallelogram with vertices ${\bf v}_1$, ${\bf v}_2$, ${\bf v}_5$ and ${\bf v}_6$, it can be deduced by convexity that
$${\rm vol}(P_8)>{\rm vol}(P)\ge 6,\eqno(56)$$
which contradicts the assumption of (53). Thus, to prove the theorem it is sufficient to deal with the three cases
$$y_2=y_1=-{3\over 2}, \ -2,\ -{5\over 2}.\eqno(57)$$

\noindent
{\bf Case 1.} $y_2=y_1=-{3\over 2}.$  In this case,
$$y_{i+1}-y_i=1$$
must hold for all $i=2, 3$ and $4$. Then, it follows by Lemma 1 that all the midpoints of $G_2$, $G_3$ and $G_4$ belong to ${1\over 2}\Lambda $. Furthermore, by a uni-modular transformation
$$\left\{
\begin{array}{ll}
x'\hspace{-0.2cm}&=x-ky,\\
y'\hspace{-0.2cm}&=y,
\end{array}\right.$$
with a suitable integer $k$, we may assume that $-{5\over 4}\le x_1<{1\over 4}$.

If $G_2$ is vertical, then $x_2$ is an integer or an half integer. Consequently, we have $x_1\in {1\over 2}\mathbb{Z}$. Therefore $x_1$ only can be $-1$, $-{1\over 2}$ or $0$. By considering three subcases with respect to $x_1=-1$, $-{1\over 2}$ or $0$, it can be deduced that there is no octagon of this type satisfying Lemma 1. For example, when $x_1=-{1\over 2}$, by Lemma 1 and convexity we have ${\bf v}_1=\left(-{1\over 2}, -{3\over 2}\right)$, ${\bf v}_2=\left({1\over 2}, -{3\over 2}\right)$, ${\bf v}_3=\left({1\over 2}, -{1\over 2}\right)$, ${\bf v}_4=\left({1\over 2}, {1\over 2}\right)$, ${\bf v}_5=-{\bf v}_1$, ${\bf v}_6=-{\bf v}_2$, ${\bf v}_7=-{\bf v}_3$ and ${\bf v}_8=-{\bf v}_4$. Then, $P_8$ is no longer an octagon but a parallelogram.

\begin{figure}[!ht]
\centering
\includegraphics[scale=0.5]{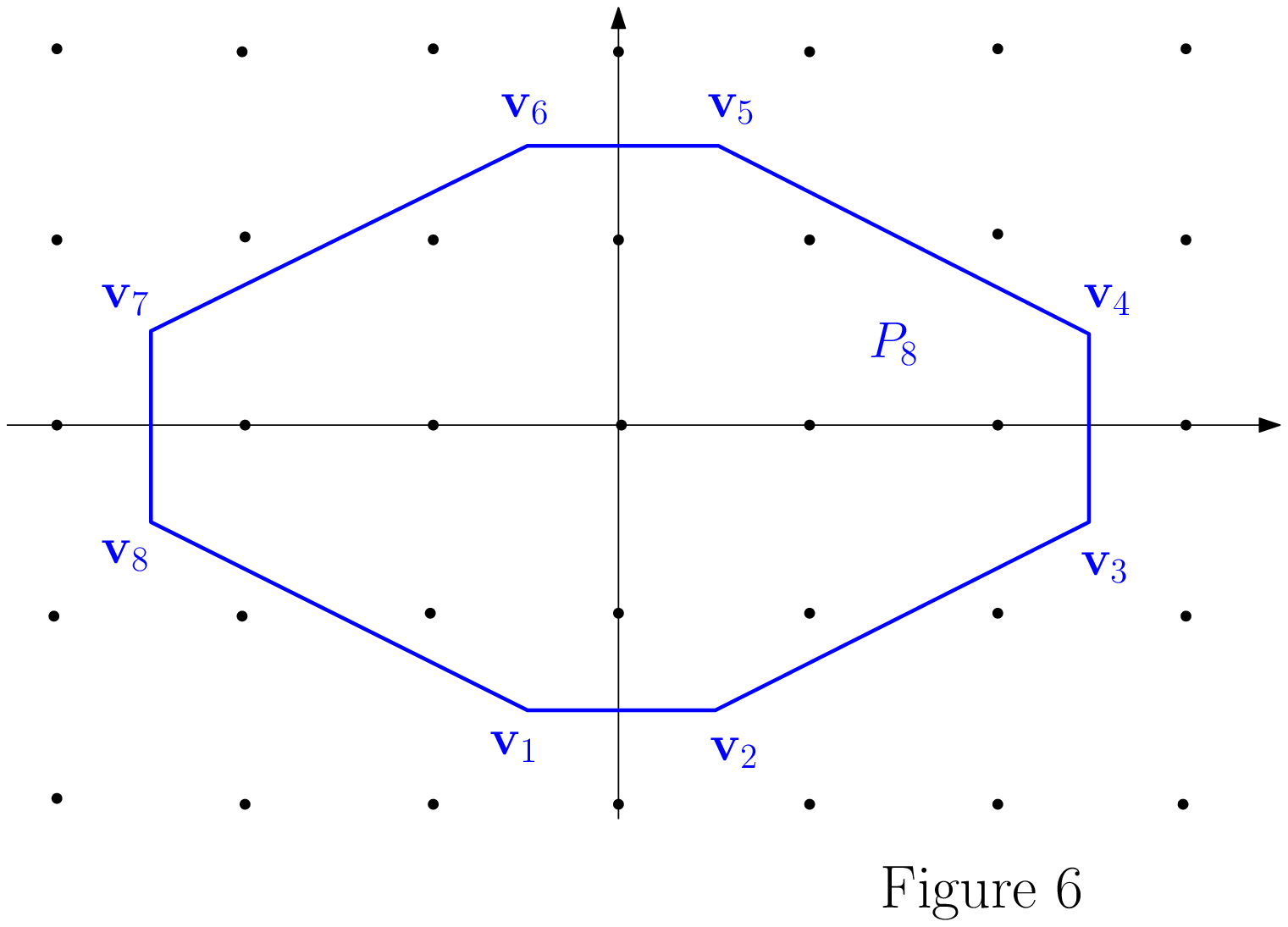}
\end{figure}

If $G_3$ is vertical, then $x_1$ must be an integer or an half integer as well. Therefore, it only can be $-1$, $-{1\over 2}$ or $0$. By considering three subcases with respect to $x_1=-1$, $-{1\over 2}$ or $0$, it can be deduced that
$${\rm vol}(P_8)\ge 7,\eqno(58)$$
which contradicts the assumption of (53). For example, when $x_1=-{1\over 2}$, by Lemma 1 and convexity we have ${\bf v}_1=\left(-{1\over 2}, -{3\over 2}\right)$, ${\bf v}_2=\left({1\over 2}, -{3\over 2}\right)$, ${\bf v}_3=\left({1\over 2}+k, -{1\over 2}\right)$, ${\bf v}_4=\left({1\over 2}+k, {1\over 2}\right)$, ${\bf v}_5=-{\bf v}_1$, ${\bf v}_6=-{\bf v}_2$, ${\bf v}_7=-{\bf v}_3$ and ${\bf v}_8=-{\bf v}_4$, where $k$ is a positive integer. Then, as shown by Figure 6, it can be deduced that
$${\rm vol}(P_8)=3+4 k\ge 7.\eqno(59)$$

If none of the three edges $G_2$, $G_3$ and $G_4$ is vertical, by convexity it is sufficient to deal with the following three subcases.

\medskip
\noindent
{\bf Subcase 1.1.}  $x'_3>\max\{ x'_2, x'_4\}$. Then we replace the eight vertices ${\bf v}_3$, ${\bf v}_4$, ${\bf v}_5$, ${\bf v}_6$, ${\bf v}_7$, ${\bf v}_8$, ${\bf v}_1$ and ${\bf v}_2$ by ${\bf v}'_3=(x'_3, -{1\over 2})$, ${\bf v}'_4=(x'_3, {1\over 2})$,
${\bf v}'_5=(2x'_4-x'_3, {3\over 2})$, ${\bf v}'_6=(2x'_4-x'_3-1, {3\over 2})$, ${\bf v}'_7=-{\bf v}'_3$, ${\bf v}'_8=-{\bf v}'_4$,  ${\bf v}'_1=-{\bf v}'_5$ and ${\bf v}'_2=-{\bf v}'_6$, respectively (as shown by Figure 7). In practice, one first makes $G_3$ vertical and then changes the other vertices successively. Clearly, this process does not change the area of the polygon. Then one can deduce that $x'_3\ge {3\over 2}$ and therefore
$${\rm vol}(P_8)=3\cdot 2x'_3-(2x'_3-1)=4x'_3+1\ge 7,\eqno (60)$$
which contradicts the assumption of (53).

\begin{figure}[!ht]
\centering
\includegraphics[scale=0.55]{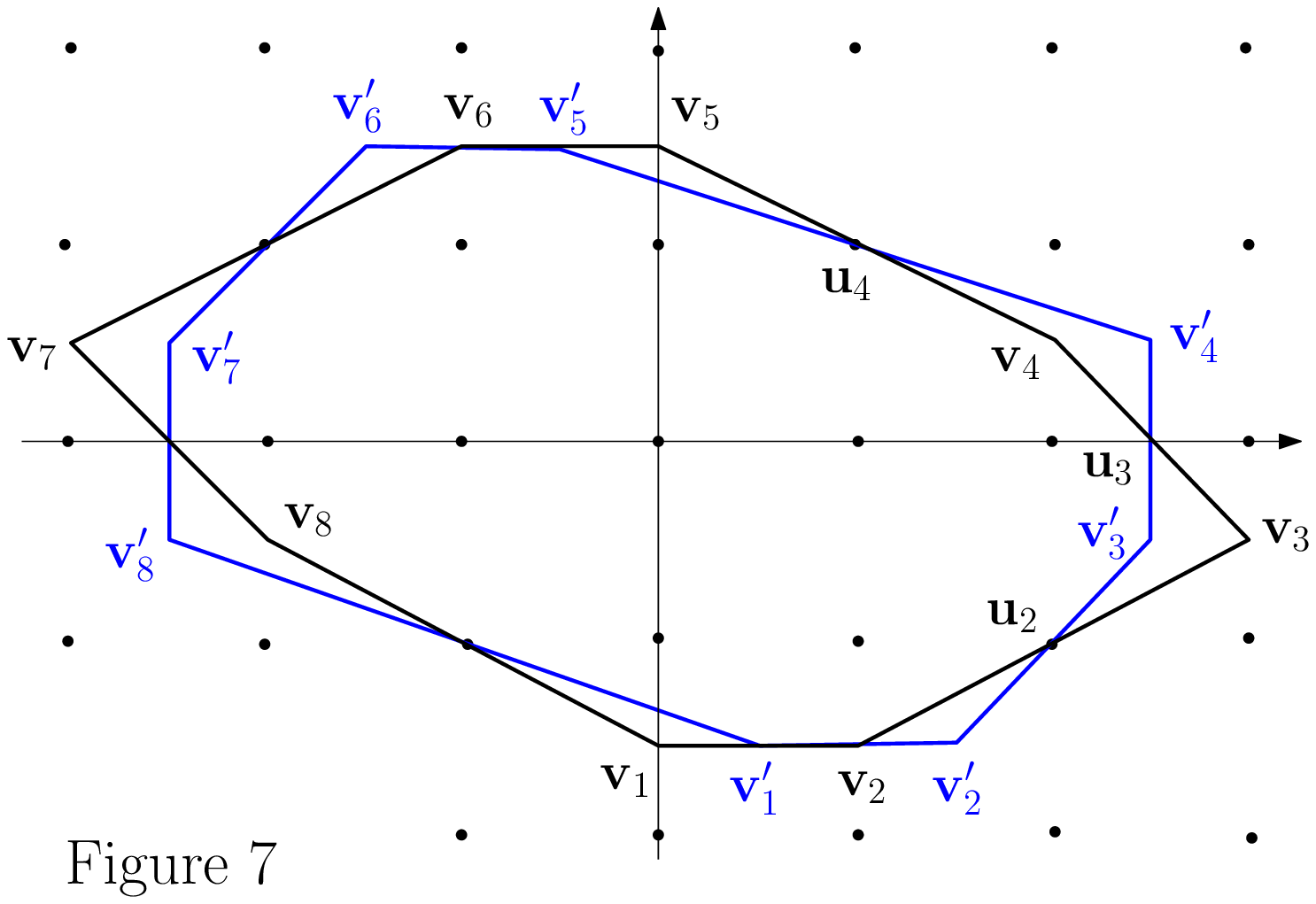}
\end{figure}

\medskip
\noindent
{\bf Subcase 1.2.} $x'_2>\max \{ x'_3, x'_4\}$. If $x_3>x_2$, one can repeat the above process. At the end we get $x'_2\ge 2$ and
$${\rm vol}(P_8)> 3\cdot 2x'_2-2(2x'_2-1)=2x'_2+2\ge 6,\eqno(61)$$
which contradicts the assumption of (53). If $x_2>x_3$, since $-{5\over 4}\le x_1<{1\over 4}$, ${\bf u}_2$ only can be $(1,-1)$, $({1\over 2}, -1)$, $(0,-1)$ or $(-{1\over 2},-1)$. Then it can be easily checked that there is no convex octagon of this type satisfying Lemma 1.

\medskip
\noindent
{\bf Subcase 1.3.} $x'_2=x'_3>x'_4$. Then, we replace the eight vertices ${\bf v}_2$, ${\bf v}_3$, ${\bf v}_4$, ${\bf v}_5$, ${\bf v}_6$, ${\bf v}_7$, ${\bf v}_8$ and  ${\bf v}_1$ by ${\bf v}'_2=(x'_2, -{3\over 2})$, ${\bf v}'_3=(x'_2, -{1\over 2})$, ${\bf v}'_4=(x'_2, {1\over 2})$,
${\bf v}'_5=2{\bf u}_4-{\bf v}'_4$, ${\bf v}'_6=-{\bf v}'_2$, ${\bf v}'_7=-{\bf v}'_3$, ${\bf v}'_8=-{\bf v}'_4$ and  ${\bf v}'_1=-{\bf v}'_5$, respectively (as shown by Figure 8). In practice, one first makes $G_2$ and $G_3$ vertical and then changes the other vertices successively, keeping the rules of Lemma 1. Clearly, this process does not change the area of the polygon, $x'_2\ge 1$ and therefore
$${\rm vol}(P_8)=3\cdot 2x'_2-(2x'_2-1)=4x'_2+1\not= 6.\eqno (62)$$

\begin{figure}[!ht]
\centering
\includegraphics[scale=0.55]{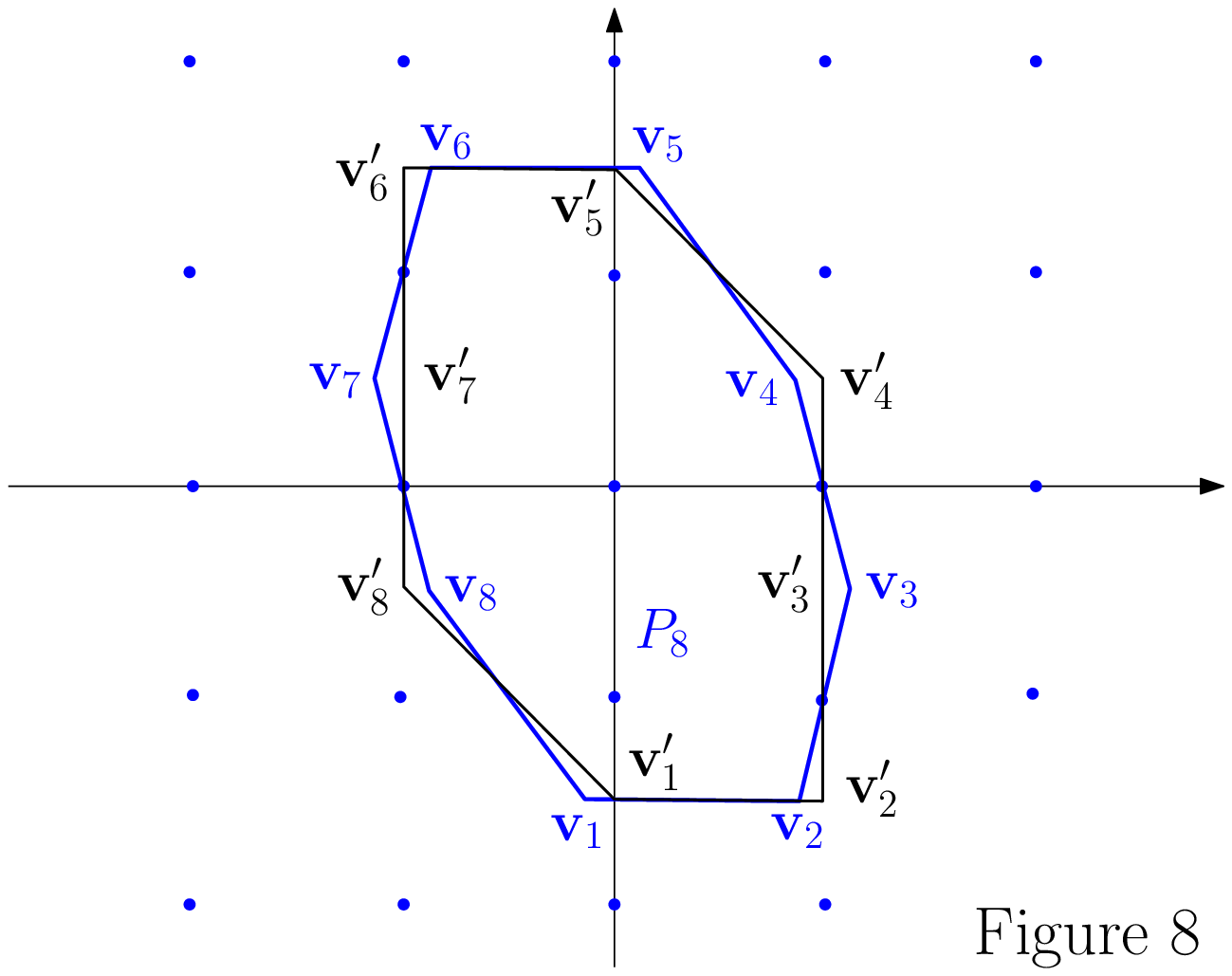}
\end{figure}

\medskip\noindent
{\bf Case 2.} $y_2=y_1=-2.$  Then, it can be deduced that one of $y_3-y_2$, $y_4-y_3$ and $y_5-y_4$ is two and the others are ones, and the midpoint ${\bf u}_i$ must belong to ${1\over 2}\Lambda $ whenever $y_{i+1}-y_i=1$. Furthermore, we may assume that $-{3\over 2}\le x_1<{1\over 2}$ by a uni-modular transformation and assume that $G_i$ is primitive if it is a lattice vector by reduction.

If one of $G_2$, $G_3$ and $G_4$ is vertical, it can be easily deduced that
$${\rm vol}(P_8)\ge 7.\eqno (63)$$
For instance, when $G_3$ is vertical, we have $x_3-x_2\ge 1$, $x_4-x_5\ge 1$ and thus $x_3=x'_3=x_4\ge {3\over 2}$. Then, it can be deduced that
$${\rm vol}(P_8)\ge 4\cdot 2x_3-2(2x_3-1)=4x_3+2\ge 8,\eqno(64)$$
which contradicts the assumption of (53).

Now, we assume that all $G_2$, $G_3$ and $G_4$ are not vertical.

\smallskip\noindent
{\bf Subcase 2.1.} {\it $y_3-y_2=2$ and ${\bf u}_2\not\in {1\over 2}\Lambda$}. Then ${\bf v}_3-{\bf v}_2=(k,2)$ is a lattice vector, where $k$ is a positive integer (When $k$ is negative, one can easily deduce that $P_8$ can not be a convex octagon). On the other hand, it follows by the assumption $-{3\over 2}\le x_1<{1\over 2}$ that
$${\bf v}_5-{\bf v}_2=(x, 4),$$
where $-2<x\le 2$. Let $P$ denote the parallelogram with vertices ${\bf v}_1$, ${\bf v}_2$, ${\bf v}_5$ and ${\bf v}_6$, and let $T$ denote the triangle with vertices ${\bf v}_2$, ${\bf v}_3$ and ${\bf v}_5$, as shown by Figure 9.

\begin{figure}[!ht]
\centering
\includegraphics[scale=0.5]{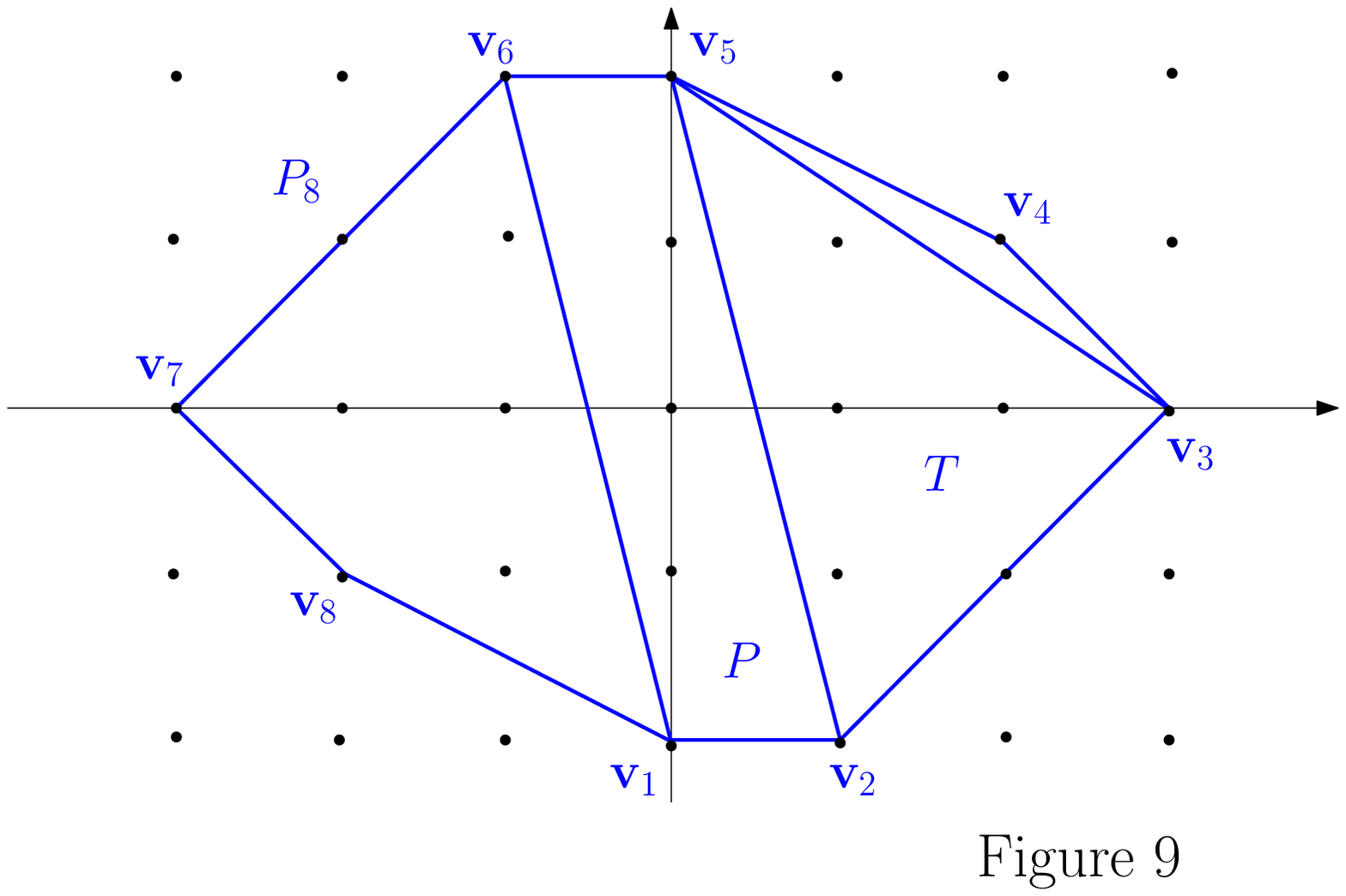}
\end{figure}

If $k\ge 2$, one can deduce
$${\rm vol}(T)={1\over 2}\left| \begin{array}{cc}
k&2\\
x&4
\end{array}\right|=2k-x\ge 2\eqno(65)$$
and therefore
$${\rm vol}(P_8)> {\rm vol}(P)+2\cdot {\rm vol}(T)\ge 8,\eqno(66)$$
which contradicts the assumption of (53).

If $k=x_3-x_2=1$, $G_2\cap {1\over 2}\Lambda \not=\emptyset $ and ${\bf u}_2\not\in {1\over 2}\Lambda $, one can deduce that $x_2\in {1\over 4}\mathbb{Z}$ and therefore $x_1\in {1\over 4}\mathbb{Z}$. In fact, by checking all the eight cases $x_1=-{3\over 2},$ $-{5\over 4},$ $-1,$ $-{3\over 4},$ $-{1\over 2},$ $-{1\over 4},$ $0$ or ${1\over 4}$, it can be shown that there is no such octagon satisfying the conditions of Lemma 1.
For example, when $x_1={1\over 4}$, by convexity (as shown by Figure 10) the only candidate for ${\bf u}_3$ is ${\bf u}'_3=(2, {1\over 2})$ and the only candidates for ${\bf u}_4$ are ${\bf u}'_4=({1\over 2}, {3\over 2})$ and ${\bf u}^*_4=(1, {3\over 2})$. However, no octagon $P_8$ satisfying Lemma 1 can be constructed from these candidate midpoints.

\begin{figure}[!ht]
\centering
\includegraphics[scale=0.5]{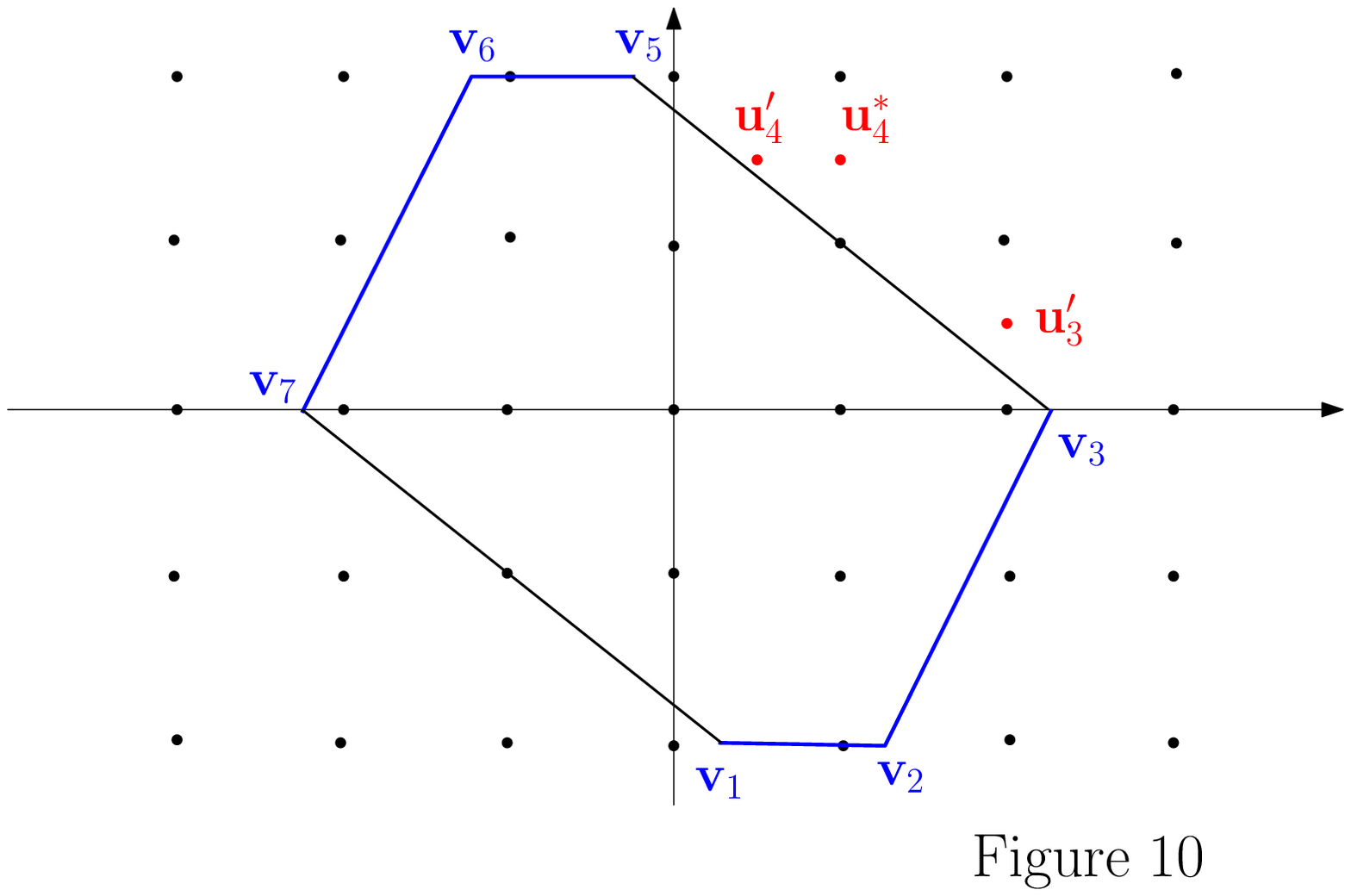}
\end{figure}

\noindent
{\bf Subcase 2.2.} {\it $y_4-y_3=2$ and ${\bf u}_3\not\in {1\over 2}\Lambda$.}  Then ${\bf v}_4-{\bf v}_3=(k,2)$ is a lattice vector, where $k$ is a positive integer (if it is negative, then make a reflection with respect to the $x$-axis). On the other hand, it follows by the assumption $-{3\over 2}\le x_1<{1\over 2}$ that
$${\bf v}_5-{\bf v}_2=(x, 4),$$
where $-2<x\le 2$. Let $P$ denote the parallelogram with vertices ${\bf v}_1$, ${\bf v}_2$, ${\bf v}_5$ and ${\bf v}_6$, and let $T$ denote the triangle with vertices ${\bf v}_2$, ${\bf v}'_3={\bf v}_2+({\bf v}_4-{\bf v}_3)$ and ${\bf v}_5$, as shown by Figure 11.

\begin{figure}[!ht]
\centering
\includegraphics[scale=0.5]{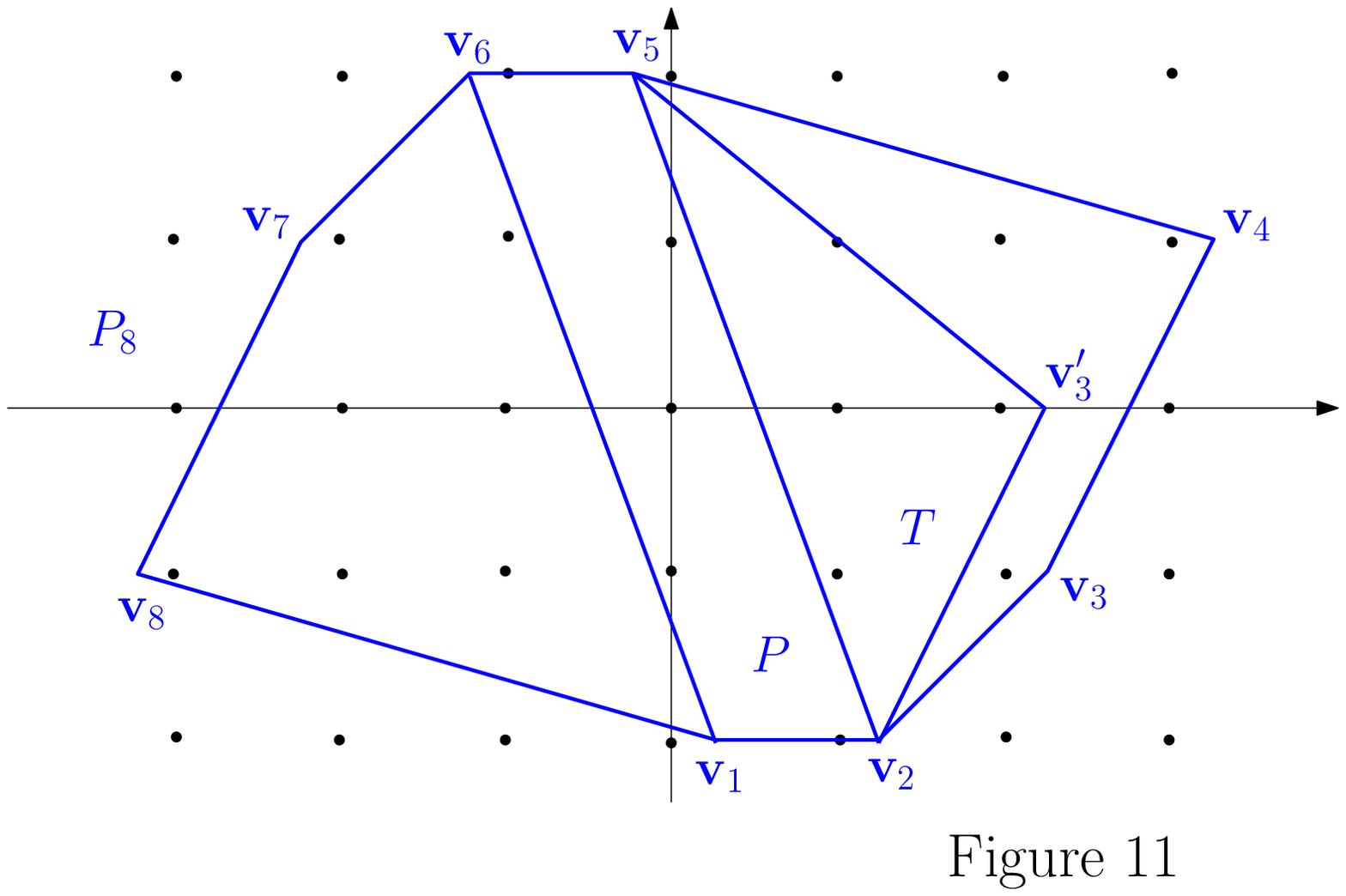}
\end{figure}

If $k\ge 2$, one can deduce
$${\rm vol}(T)={1\over 2}\left| \begin{array}{cc}
k&2\\
x&4
\end{array}\right|=2k-x\ge 2\eqno(67)$$
and therefore
$${\rm vol}(P_8)> {\rm vol}(P)+2\cdot {\rm vol}(T)\ge 8,\eqno(68)$$
which contradicts the assumption of (53).

If $k=x_4-x_3=1$, $G_3\cap {1\over 2}\Lambda \not=\emptyset $ and ${\bf u}_3\not\in {1\over 2}\Lambda $, one can deduce that $x_3\in {1\over 4}\mathbb{Z}$ and therefore $x_1\in {1\over 4}\mathbb{Z}$. By checking all the eight cases $x_1=-{3\over 2},$ $-{5\over 4},$ $-1,$ $-{3\over 4},$ $-{1\over 2},$ $-{1\over 4},$ $0$ or ${1\over 4}$, it can be deduced that
$${\rm vol}(P_8)\ge 7.\eqno (69)$$
For example, when $x_1=-{3\over 2}$, we define ${\bf v}'_3=({3\over 2}, -1)$, ${\bf v}'_4=({5\over 2}, 1)$, ${\bf v}'_7=(-{3\over 2},1)$,
${\bf v}'_8=(-{5\over 2},-1)$, and define $P'_8$ to be the octagon with vertices ${\bf v}_1$, ${\bf v}_2$, ${\bf v}'_3$, ${\bf v}'_4$, ${\bf v}_5$, ${\bf v}_6$, ${\bf v}'_7$ and ${\bf v}'_8$, as shown by Figure 12. By shifting $G_3$ and $G_7$, one can deduce $P'_8\subseteq P_8$ and therefore
$${\rm vol}(P_8)\ge {\rm vol}(P'_8)=13.\eqno(70)$$

\begin{figure}[!ht]
\centering
\includegraphics[scale=0.5]{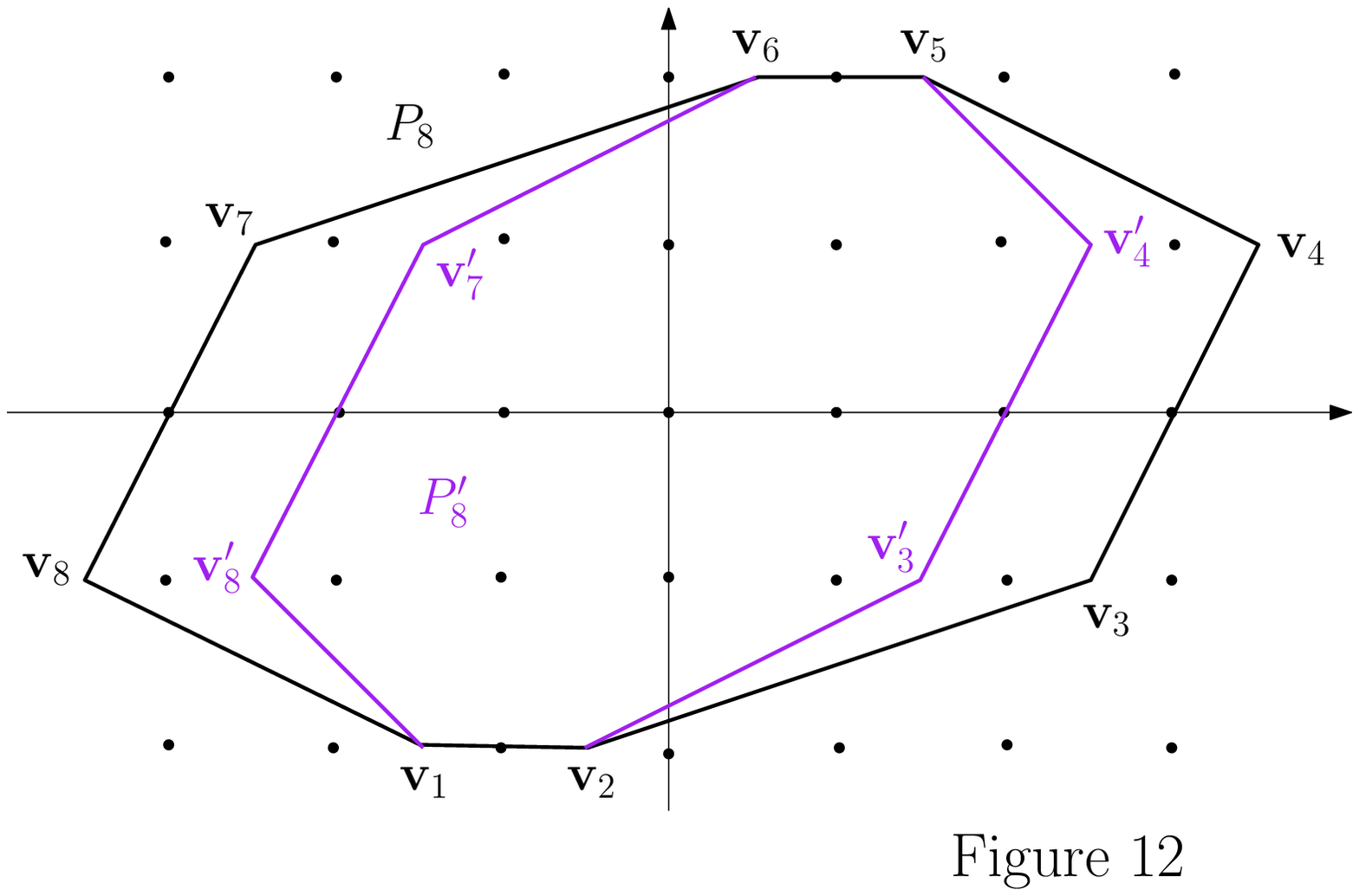}
\end{figure}

\noindent
{\bf Subcase 2.3.} {\it None of the three edges $G_2$, $G_3$ and $G_4$ is vertical and all ${\bf u}_2$, ${\bf u}_3$ and ${\bf u}_4$ belong to ${1\over 2}\Lambda $.} Then, it is sufficient to consider the following three situations.

\medskip\noindent
{\bf Subcase 2.3.1.}  $x'_3>\max\{ x'_2, x'_4\}$. Similar to Subcase 1.1, we get $x'_3\ge {3\over 2}$ and therefore
$${\rm vol}(P_8)\ge 4\cdot 2x'_3-2(2x'_3-1)=4x'_3+2\ge 8,\eqno (71)$$
which contradicts the assumption of (53).

\medskip
\noindent
{\bf Subcase 2.3.2.} $x'_2>\max \{ x'_3, x'_4\}$. If $x_3>x_2$, just like Subcase 1.2, one can get $x'_2\ge {3\over 2}$ and
$${\rm vol}(P_8)> 4\cdot 2x'_2-3(2x'_2-1)\ge 6,\eqno(72)$$
which contradicts the assumption of (53).

If $x_2>x_3$ and $y_3-y_2=1$, since $-{3\over 2}\le x_1<{1\over 2}$, ${\bf u}_2$ only can be $(1, -{3\over 2})$, $({1\over 2}, -{3\over 2})$,
$(0, -{3\over 2})$ or $(-{1\over 2}, -{3\over 2})$. Then it can be routinely checked that there is no convex octagon of this type satisfying Lemma 1.

\begin{figure}[!ht]
\centering
\includegraphics[scale=0.5]{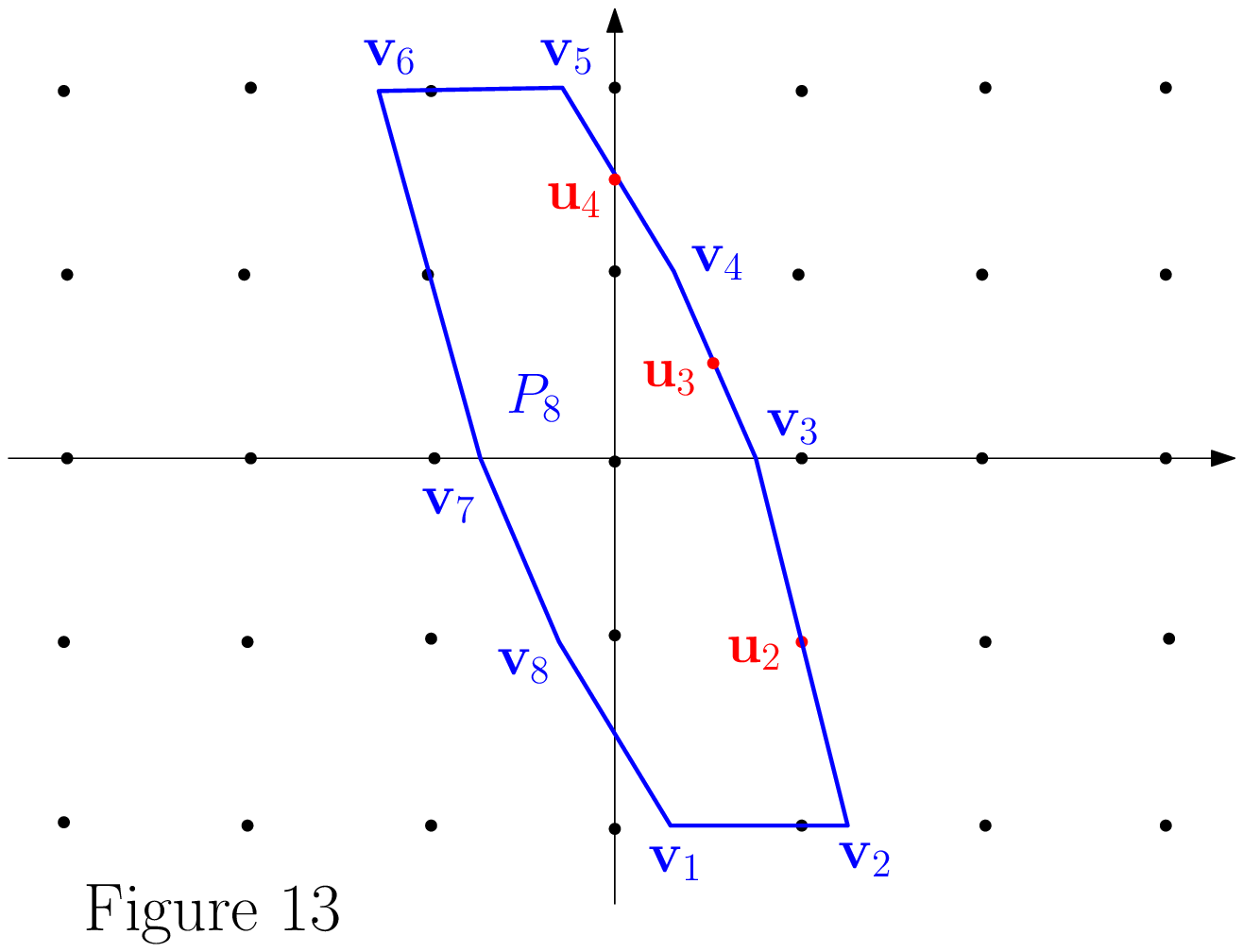}
\end{figure}

If $x_2>x_3$ and $y_3-y_2=2$, since $-{3\over 2}\le x_1<{1\over 2}$, ${\bf u}_2$ only can be $(1, -1)$, $({1\over 2}, -1)$,
$(0, -1)$ or $(-{1\over 2}, -1)$. By checking these four cases, it can be shown that there is only one class of such convex octagons satisfying Lemma 1. Namely, the ones satisfying ${\bf u}_2=(1,-1)$, ${\bf u}_3=({1\over 2}, {1\over 2})$ and ${\bf u}_4=(0,{3\over 2})$, as shown in Figure 13. In other words, they are the octagons with vertices ${\bf v}_1=(\beta , -2)$, ${\bf v}_2=(1+\beta , -2)$, ${\bf v}_3=(1-\beta , 0)$, ${\bf v}_4=(\beta , 1)$, ${\bf v}_5=-{\bf v}_1$, ${\bf v}_6=-{\bf v}_2$, ${\bf v}_7=-{\bf v}_3$, ${\bf v}_8=-{\bf v}_4$, where ${1\over 4}<\beta <{1\over 3}$. Then, one can deduce that
$${\rm vol}(P_8)=5.\eqno(73)$$

\medskip
\noindent
{\bf Subcase 2.3.3.} $x'_2=x'_3>x'_4$. Similar to Subcase 1.3, one can deduce $x'_2\ge 1$ and therefore
$${\rm vol}(P_8)\ge 4\cdot 2x'_3-2(2x'_3-1)=4x'_3+2\ge 6,\eqno (74)$$
where the equalities hold if and only if $P_8$ with vertices ${\bf v}_1=(\alpha-1 , 2)$, ${\bf v}_2=(\alpha , -2)$, ${\bf v}_3=(1-\alpha , 0)$, ${\bf v}_4=(1+\alpha , -1)$, ${\bf v}_5=-{\bf v}_1$, ${\bf v}_6=-{\bf v}_2$, ${\bf v}_7=-{\bf v}_3$ and ${\bf v}_8=-{\bf v}_4$, where $0<\alpha <{1\over 6}$ (as shown in Figure 14).

\begin{figure}[!ht]
\centering
\includegraphics[scale=0.5]{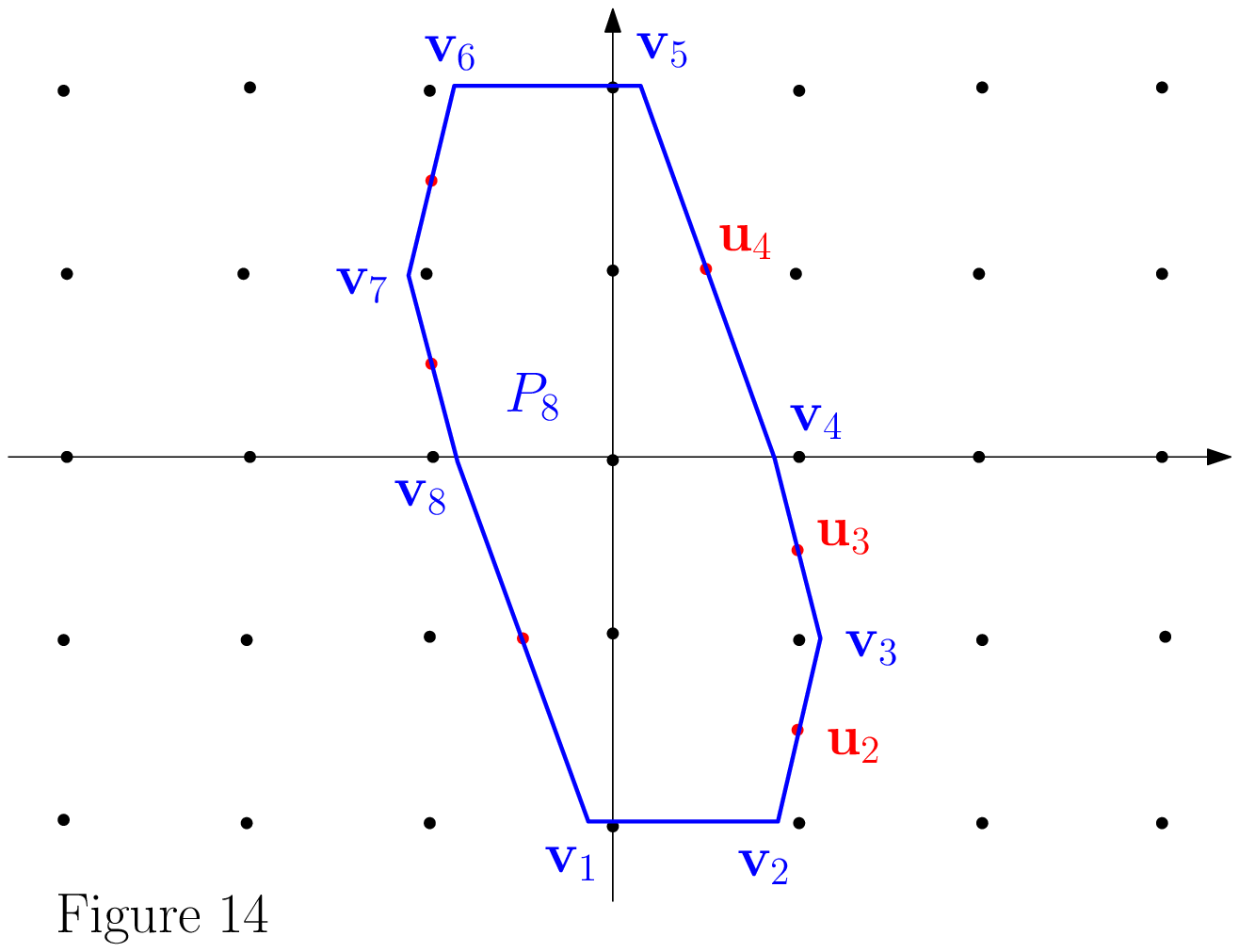}
\end{figure}

\medskip
\noindent
{\bf Case 3.} $y'_1= -{5\over 2}$. Then all $y_i-y_{i+1}$ are positive integers for $2\le i\le 4$. If
$${\bf u}_i\in \mbox{${1\over 2}$} \Lambda $$
hold for all $i=2,$ $3$ and $4$, similar to Subcase 1.1, Subcase 1.2 and Subcase 1.3 one can deduce
$$\tau^*(P_8)={\rm vol}(P_8)\ge 7.\eqno(75)$$

If ${\bf u}_i\not\in {1\over 2}\Lambda$ holds for one of these indices, then we have $y_i-y_{i+1}=2$ or $3$. By a uni-modular transformation, we may assume that $-{7\over 4}\le x_1<{3\over 4}.$ Then we have ${\bf v}_2-{\bf v}_6=(x, 5)$, where $-{5\over 2}\le x< {5\over 2}$. If ${\bf v}_i-{\bf v}_{i+1}=(k,2)$ with $|k|\ge 2$, let $Q$ denote the pentagon with vertices ${\bf v}_2$, ${\bf v}_3$, ${\bf v}_4$, ${\bf v}_5$ and ${\bf v}_6$, then we have
$${\rm vol}(Q)> {1\over 2}\left|
\begin{array}{cc}
x & 5\\
k& 2
\end{array}\right| = {1\over 2}\left| 2x-5k\right|\ge {5\over 2}\eqno(76)$$
and thus
$$\tau^*(P_8)={\rm vol}(P)+2\cdot {\rm vol}(Q)\ge 10.\eqno(77)$$
If ${\bf v}_i-{\bf v}_{i+1}=(k,2)$ with $k=\pm 1$, then we have $x_1\in {1\over 4}\mathbb{Z}$ and therefore $x\in {1\over 2}\mathbb{Z}$ and $-{5\over 2}\le x\le 2$. By considering two subcases with respect to $x_1=-{7\over 4}$ and $x_1\not=-{7\over 4}$, we can get
$${\rm vol} (Q)>{1\over 2}\eqno(78)$$
and
$$\tau^*(P_8)={\rm vol}(P)+2\cdot {\rm vol}(Q)>6.\eqno(79)$$
The $y_i-y_{i+1}=3$ case can be eliminated in a similar way.

As a conclusion of all these cases, Lemma 10 is proved. \hfill{$\Box$}

\vspace{0.6cm}
\noindent
{\Large\bf 4. Proofs of the Theorems}

\bigskip\noindent
{\bf Proof of Theorem 1.} Theorem 1 follows from Lemmas 5, 7-10 immediately. \hfill{$\Box$}

\bigskip\noindent
{\bf Proof of Theorem 2.} Follow the proof of Theorem 2 in Zong \cite{zong}, this theorem can be easily proved. \hfill{$\Box$}

\vspace{0.6cm}\noindent
{\bf Acknowledgements.} This work is supported by 973 Program 2013CB834201.

\bibliographystyle{amsplain}

\begin{thebibliography}{99}
\bibitem{alek}A. D. Aleksandrov, On tiling space by polytopes, {\it Vestnik Leningrad Univ. Ser. Mat. Fiz. Him}. {\bf 9} (1954), 33-43.
\bibitem{boll}U. Bolle, On multiple tiles in $R^2$, {\it Intuitive Geometry,} Colloq. Math. Soc. J. Bolyai {\bf 63}, North-Holland, Amsterdam, 1994.
\bibitem{delo}B. N. Delone, Sur la partition reguli$\grave{e}$re de l'espace $\grave{a}$ $4$ dimensions I, II,
{\it Izv. Akad. Nauk SSSR, Ser. VII} (1929), 79-110; 147-164.
\bibitem{dgsw}M. Dutour Sikiri$\acute{\rm c}$, A. Garber, A. Sch$\ddot{\rm u}$rmann and C. Waldmann, The complete classification of five-dimensional Dirichlet-Voronoi polyhedra of translational lattices, {\it Acta Crystallographica} {\bf A72} (2016), 673-683.
\bibitem{enge}P. Engel, On the symmetry classification of the four-dimensional parallelohedra, {\it Z. Kristallographie} {\bf 200} (1992), 199-213.
\bibitem{fedo}E. S. Fedorov, Elements of the study of figures, {\it Zap. Mineral. Imper. S. Petersburgskogo Ob$\check{s}$$\check{c}$},
{\bf 21}(2) (1885), 1-279.
\bibitem{furt}P. Furtw\"angler, \"Uber Gitter konstanter Dichte, {\it Monatsh. Math. Phys.} {\bf 43} (1936), 281-288.
\bibitem{grs}N. Gravin, S. Robins and D. Shiryaev, Translational tilings by a polytope, with multiplicity. {\it Combinatorica}
{\bf 32} (2012), 629-649.
\bibitem{gkrs}N. Gravin, M. N. Kolountzakis, S. Robins and D. Shiryaev, Structure results for multiple tilings in 3D.
{\it Discrete Comput. Geom.} {\bf 50} (2013), 1033-1050.
\bibitem{grub}P. M. Gruber and C. G. Lekkerkerker, {\it Geometry of Numbers} (2nd ed.), North-Holland, 1987.
\bibitem{hajo}G. Haj\'os, \"Uber einfache und mehrfache Bedeckung des $n$-dimensionalen Raumes mit einem W\"urfelgitter, {\it Math. Z.}
{\bf 47} (1941), 427-467.
\bibitem{kolo}M. N. Kolountzakis, On the structure of multiple translational tilings by polygonal regions. {\it Discrete Comput. Geom}.
{\bf 23} (2000), 537-553.
\bibitem{lazo}J. C. Lagarias and C. Zong, Mysteries in packing regular tetrahedra, {\it Notices Amer. Math. Soc.}
{\bf 59} (2012), 1540-1549.
\bibitem{mann}C. Mann, J. McLoud-Mann and D. Von Derau, Convex pentagons that admit $i$-block transitive tilings, {\it Geom. Dedicata} {\bf 194} (2018), 141-167.
\bibitem{mcmu}P. McMullen, Convex bodies which tiles space by translation, {\it Mathematika} {\bf 27} (1980), 113-121.
\bibitem{mink}H. Minkowski, Allgemeine Lehrs\"atze \"uber konvexen Polyeder, {\it Nachr. K. Ges. Wiss. G\"ottingen, Math.-Phys. KL}. (1897), 198-219.
\bibitem{rao}M. Rao, Exhaustive search of convex pentagons which tile the plane, arXiv:1708.00274.
\bibitem{rein}K. Reinhardt, \"Uber die Zerlegung der Ebene in Polygone, {\it Dissertation,} Universit\"at Frankfurt am Main, 1918.
\bibitem{robi}R. M. Robinson, Multiple tilings of $n$-dimensional space by unit cubes, {\it Math. Z.} {\bf 166} (1979), 225-275.
\bibitem{stog}M. I. $\check{S}$togrin, Regular Dirichlet-Voronoi partitions for the second triclinic group (in Russian),
{\it Proc. Steklov. Inst. Math.} {\bf 123} (1975).
\bibitem{venk}B. A. Venkov, On a class of Euclidean polytopes, {\it Vestnik Leningrad Univ. Ser. Mat. Fiz. Him}. {\bf 9} (1954), 11-31.
\bibitem{voro}G. F. Voronoi,  Nouvelles applications des paramm\`etres continus \`a la th\'eorie des formes quadratiques.
Deuxi\`eme M\'emoire. Recherches sur les parall\'elo\`edres primitifs, {\it J. reine angew. Math.} {\bf 134} (1908), 198-287; {\bf 135} (1909), 67-181.
\bibitem{yz1}Q. Yang and C. Zong, Multiple lattice tilings in Euclidean spaces, {\it Canadian Math. Bull.} in press.
\bibitem{yz2}Q. Yang and C. Zong, Characterization of the two-dimensional five-fold translative tiles, arXiv:1711.02514.
\bibitem{zong05}C. Zong, What is known about unit cubes. {\it Bull. Amer. Math. Soc.} {\bf 42} (2005), 181-211.
\bibitem{zong06}C. Zong, {\it The Cube: A Window to Convex and Discrete Geometry.} Cambridge University Press, Cambridge, 2006.
\bibitem{zong14}C. Zong, Packing, covering and tiling in two-dimensional spaces, {\it Expo. Math.} {\bf 32} (2014), 297-364.
\bibitem{zong}C. Zong, Characterization of the two-dimensional five-fold lattice tiles, arXiv:1712.01122.

\end{thebibliography}

\vspace{0.6cm}
\noindent
Chuanming Zong, Center for Applied Mathematics, Tianjin University, Tianjin 300072, China

\noindent
Email: cmzong@math.pku.edu.cn

\end{document}